\documentclass[11pt]{article}
\usepackage{amssymb,amsfonts,amsmath,amsthm,mathrsfs,enumerate,enumitem}

\newtheorem{Theorem}{Theorem}[section]
\newtheorem{Lemma}[Theorem]{Lemma}
\newtheorem{Corollary}[Theorem]{Corollary}
\newtheorem{prop}[Theorem]{Proposition}
\newtheorem{rem}[Theorem]{Remark}

\newtheorem{hyp}{Assumption}[section]

\usepackage{color}
\usepackage[colorlinks = true,
linkcolor = blue,
urlcolor  = blue,
citecolor = blue,
anchorcolor = blue]{hyperref}
\usepackage{graphicx}

\numberwithin{equation}{section}

\usepackage{empheq}
\DeclareMathOperator{\re}{Re\,}
\title{\bf A new perspective of exponential stability for Timoshenko systems under history and thermal effects\thanks{This paper is dedicated to the   memory of Pinheiro's younger sister, Gabriela B. Pinheiro, who had already decided to study mathematics but whose life was taken at age 14 by  early cancer.}}
\author{
	\small {\bf Marcio A. Jorge Silva}
	\thanks{Corresponding author. Email: \href{marcioajs@uel.br}{ marcioajs@uel.br} }
	\thanks{Partially  supported by the CNPq, Grant \#301116/2019-9.}\\
	\small Department of Mathematics, State University of Londrina,\\
	\small  Londrina 86057-970, Paran\'a, Brazil. \medskip
	\\
	\small {\bf Sandro B. Pinheiro} \thanks{Supported by the CAPES, Finance Code 001. 	(Master and Ph.D. Scholarships)
	}\\
	\small Department of Mathematics, State University of Maring\'a,\\
	\small Maring\'a 87020-900, Paran\'a,  Brazil.
}

\date{ }

\begin{document}

\maketitle

\begin{abstract}
We address a  Timoshenko system  with memory in the history context and thermoelasticity of	
type III  for heat conduction.  Our main goal is to prove its uniform (exponential) stability  by illustrating carefully the sensitivity of the heat and history couplings on the  Timoshenko system. This investigation
 contrasts previous insights on the subject and promotes a new perspective with respect to the stability of the thermo-viscoelastic problem carried out, by combining the whole strength of   history and thermal effects.
\end{abstract}
 
\noindent{\bf Keywords:} Timoshenko system; Exponential stability; History; Type III thermoelasticity.
 
\smallskip

\noindent{\bf 2020 MSC:}  35B35; 35B40; 35Q79; 74D05; 74F05; 74H40.

%   
%\begin{quote}
%{\small	\tableofcontents}
%\end{quote}

\section{Introduction}\label{sec-introdu}

In the present article, our main goal is to study the uniform stability of the following Timoshenko beam model with  
thermoelasticity of type III and memory with history 
\begin{align}
	\rho_1\varphi_{tt}  - k(\varphi_{x}+\psi)_{x} + \sigma\theta_{tx}\, = \, 0 & \  \ \mbox{ in} \  \ (0,l)\times(0,\infty), \label{intro 1.5}\\
	\rho_2\psi_{tt}  -  \beta\psi_{xx} +k(\varphi_{x}+\psi) -\int_0^{\infty}{g(s) \eta_{xx}(s)ds} - \sigma\theta_t \, = \, 0 &\  \ \mbox{ in} \  \ (0,l)\times(0,\infty),  \label{intro 1.6} \\
	\rho_3 \theta_{tt}  - \delta\theta_{xx}  - \gamma \theta_{xxt} + \sigma (\varphi_{x}+\psi)_{t} \, = \, 0 &\  \ \mbox{ in} \  \ (0,l)\times(0,\infty), \label{intto 1.7}
	\medskip \\
	\eta_t + \eta_s -\psi_t\,=\, 0 &\  \ \mbox{ in}\  \ (0,l) \times (0,\infty)^2, \label{intro 1.8}
\end{align}
subject to initial conditions
\begin{equation} \label{init-cond}
\begin{array}{c}
(\varphi,\varphi_t,\psi,\psi_t,	\theta,\theta_t)(x,0)=(\varphi_0, \varphi_1,\psi_0, \psi_1, \theta_0, \theta_1)(x), \ \ x\in(0,l), \smallskip \\
\ \ \eta(x,0,s)= \eta_0(x,s), 
\ \eta(x,t,0)=0, \ \ x\in(0,l), \, s>0, \, t\geq0,
	\end{array}
\end{equation}
and either boundary conditions the full Dirichlet case
\begin{subequations}\label{bc}
	\begin{eqnarray}\label{bc-a}
		\left.\begin{array}{c}
			\varphi(0,t)=\varphi(l,t)=\psi(0,t)=\psi(l,t)=\theta(0,t)=\theta(l,t)  = 0, \ \ t\geq0, \smallskip \\
			\eta(0,t,s)=\eta(l,t,s)=0, \ \   t\geq 0, \ s>0,
		\end{array}\right.
	\end{eqnarray}
	or the mixed Neumann-Dirichlet one
	\begin{eqnarray}\label{bc-b}
		\left.\begin{array}{c}
			\varphi_x(0,t)=\varphi_x(l,t)=\psi(0,t)=\psi(l,t)=\theta(0,t)=\theta(l,t)  = 0, \ \ t\geq0, \smallskip \\
			\eta(0,t,s)=\eta(l,t,s)=0, \ \    t\geq 0, \ s>0, 
		\end{array}\right.
	\end{eqnarray}
\end{subequations}
where  $\rho_1, \rho_2, \rho_3, k, b, \delta, \gamma, \, \sigma,$  and $ \beta:=b-\int_{0}^{\infty}{g(s)ds} $ are positive constants whose physical meanings are very well described, the unknown functions $\varphi=\varphi(x,t), \psi=\psi(x,t), \theta=\theta(x,t),$ and  $\eta=(x,t,s)$ 
are, respectively, related to  the transversal displacement,  the rotation angle, the  temperature, and the  relative displacement history of a beam  with length $l>0$. 
The physical  details around problem \eqref{intro 1.5}-\eqref{bc} 
will be clarified in Section  \ref{sec-rise-problem}.

In order to prove the uniform stabilization of \eqref{intro 1.5}-\eqref{bc}, we consider the standard exponential assumption on the memory kernel $g$ as follows.
\begin{hyp}\label{Hyp-g}
Let us suppose that 	$g\in C^1(0,\infty)\cap L^1(0,\infty) $ satisfies
	\begin{equation}\label{sobreg}
	0<g(0)<\infty, \ 	0< \int_{0}^{\infty}{g(s)ds}< b, \ \mbox{ and } \ 0< k_1 g(s)\leq - g'(s),\; s\in (0,\infty),
	\end{equation}
	for some $k_1>0$.
\end{hyp}

Under  Assumption \ref{Hyp-g},
we shall prove that 
model \eqref{intro 1.5}-\eqref{bc} is always exponentially stable, independently of any relationship among the coefficients  and  the boundary conditions (\eqref{bc-a} or \eqref{bc-b}) taken into account (Theorem \ref{theo-expo3}).
This fact seems to be, somehow, a surprising result with respect to the possible roles of history and heat conduction in type III thermoelasticity since it strongly contrasts earlier expectations on the stability of the PDE system \eqref{intro 1.5}-\eqref{intro 1.8}. Indeed, as stated without proof in \cite[Sect. 7]{dilberto-mauro} an initial-boundary value problem related to \eqref{intro 1.5}-\eqref{bc}
{\it is  not exponentially stable} in general. The authors claimed that the stability finally depends upon the 
equal speeds of wave propagation 
\begin{equation}\label{ews}
\frac{k}{\rho_{1}} = \frac{b}{\rho_{2}},
\end{equation}
by suggesting that the model {\it has an optimal polynomial decay rate when \eqref{ews} does not hold.} More precise details on this statement will be clarified in Remark \ref{rem-santos} right after the proofs of our main result (Theorem \ref{theo-expo3}). Therefore, by means of our main result, a new perspective of stability to 
 \eqref{intro 1.5}-\eqref{bc} is provided. Surprisingly or not, we advance that the memory and temperature components produce enough dissipation in order to stabilize the whole system. In other words, one can say that problem  \eqref{intro 1.5}-\eqref{bc}  is, indeed, fully damped.
  This fact will be physically clearer 
 at the end of Section \ref{sec-rise-problem} and will be mathematically proved in Section \ref{sec-thermovisco-problem}. Under these statements, a simple question  so arises.
  \begin{itemize}
 	\item[${\bf Q1.}$]   
 	Why does our result contrast the statement given in Section 7 of \cite{dilberto-mauro}? 
 \end{itemize}

The main reason is behind the thermal and viscoelastic  couplings on the canonical Timoshenko system. In order to give a deeper response to this question, let us consider
 another 
close  Timoshenko system (but not the same)  with  
thermoelasticity of type III and memory with history 
 \begin{align}
 	\rho_1\varphi_{tt}  - k(\varphi_{x}+\psi)_{x} \, = \, 0 & \  \ \mbox{ in} \  \ (0,l)\times(0,\infty), \label{intro 1.5-n}\\
 	\rho_2\psi_{tt}  -  \beta\psi_{xx} +k(\varphi_{x}+\psi) -\int_0^{\infty}{g(s) \eta_{xx}(s)ds} + \sigma\theta_{tx} \, = \, 0 &\  \ \mbox{ in} \  \ (0,l)\times(0,\infty),  \label{intro 1.6-n} \\
 	\rho_3 \theta_{tt}  - \delta\theta_{xx}  - \gamma \theta_{xxt} + \sigma \psi_{xt} \, = \, 0 &\  \ \mbox{ in} \  \ (0,l)\times(0,\infty), \label{intto 1.7-n}
 	\medskip \\
 	\eta_t + \eta_s -\psi_t\,=\, 0 &\  \ \mbox{ in}\  \ (0,l) \times (0,\infty)^2. \label{intro 1.8-n}
 \end{align}
  
System \eqref{intro 1.5-n}-\eqref{intro 1.8-n} was addressed in \cite{salim-belkacem-09} with proper initial-boundary conditions, see    
 (2.1)-(2.4) therein. Summarizing, for exponential kernels  $g$   as in \eqref{sobreg}, the authors prove:   $(i)$ under the assumption \eqref{ews}, the energy goes to zero exponentially when $t$ goes to infinity (cf. \cite[Thm. 2.1]{salim-belkacem-09}); $(ii)$ Otherwise, if \eqref{ews} does not hold, then the energy is only semi-uniformly stable with polynomial decay-type only for regular initial data (cf. \cite[Thm. 3.1]{salim-belkacem-09}). 
 A complete characterization of the stability for \eqref{intro 1.5-n}-\eqref{intro 1.8-n} would be  provided if the authors showed lack of exponential stability in case \eqref{ews} fails, although the authors do not consider this part. In conclusion, problem \eqref{intro 1.5-n}-\eqref{intro 1.8-n} is only partially dissipative, unlike 
 \eqref{intro 1.5}-\eqref{intro 1.8}.
 Why?
 
 As a matter of fact, instead of answering the previous question ${\bf Q1}$, the  aforementioned statements drive us to another  intriguing one, as   contextualized below.  Formally speaking, both problems  \eqref{intro 1.5}-\eqref{intro 1.8} and \eqref{intro 1.5-n}-\eqref{intro 1.8-n} have the same couple of variables, the same number of equations, the same {\it  quantity} of damping  terms, the same energy, and also the same energy derivative which can be achieved with standard computations. Thus, one can ask: 
  \begin{itemize}
	\item[${\bf Q2.}$]   
	Why  is problem  
	\eqref{intro 1.5}-\eqref{intro 1.8} fully damped   while   
	\eqref{intro 1.5-n}-\eqref{intro 1.8-n}
	is partially dissipative? 
\end{itemize}
 
In order to give the precise answers to ${\bf Q1}$ and ${\bf Q2}$, we must go back to the governing equations for
Timoshenko  beams along with thermo-(visco-)elastic constitutive laws on the forces of the system. The latter are given by the bending moment and the shear force. Proceeding in this way, we can provide an explicit formulation for both problems by showing that the first one is obtained with thermal and viscoelastic  couplings on different forces of the system whereas the second one is reached by considering history and heat couplings on the same force. The full description of the latter statement is  presented in Section \ref{sec-rise-problem} and answers question  ${\bf Q2}$, at least physically. To this end, we   follow the modeling provided in \cite{michele-JEE,michele-SIAM}. Additionally, the mathematical (and technical)  answers to questions ${\bf Q1}$-${\bf Q2}$ concerning problem \eqref{intro 1.5}-\eqref{intro 1.8} are fully provided in Section \ref{sec-thermovisco-problem}. To this purpose,  we employ a refined resolvent analysis in combination with 
auxiliary results  given in Appendix \ref{sec-append-aux-results} and well-known results in linear semigroup theory.

To sum up briefly, our contributions in the present paper are: 
\begin{itemize}
	\item to bring new perspectives in what concerns  the role of memory with  history and heat conduction in type III thermoelasticity  for Timoshenko-type systems;
	
	\item to provide the correct uniform stability  result with respect to problem \eqref{intro 1.5}-\eqref{intro 1.8} subject to initial-boundary conditions \eqref{init-cond}-\eqref{bc}.
\end{itemize}

\section{Rise of the model with history and thermal effects}\label{sec-rise-problem}

Let us initially consider the 
governing equations for  
Timoshenko  beams  (cf. \cite{timoshenko,timoshenko-1}):
\begin{equation}\label{tim1}
	\left\{\begin{array}{lcl}
		\rho_1\varphi_{tt}-S_x=0 & \mbox{in} & (0,l)\times (0,\infty) , \medskip \\ 
		\rho_2\psi_{tt}-M_x+S=0  & \mbox{in} & (0,l)\times (0,\infty),
	\end{array} \right.
\end{equation}
where $\rho_{1},\rho_{2}>0$ are   physically well-known constants, $\varphi = \varphi(x,t)$,
 $\psi=\psi(x,t)$, $S=S(x,t)$, and  $M=S(x,t)$ 
 stand for the transversal displacement, rotation angle, shear force, and  bending moment, respectively. In the classical elastic case, the following constitutive laws are in place
 \begin{equation}\label{shear-bending-elastic}
 	\left\{\begin{array}{l}
 		S  \,= \,  k\big(\varphi_{x}+\psi\big) , \medskip \\ 
 		M  \,= \, b  \, \psi_{x }  ,
 	\end{array} \right.
 \end{equation}
where $k,b>0$ are again well-known constants coming from physical concerns. Therefore, replacing \eqref{shear-bending-elastic} in \eqref{tim1} we obtain the classical conservative model in differential equations for vibrations of thin beams as originated in Timoshenko's works.

Now, on the one hand, for materials containing hereditary (history) properties, the Boltzmann theory for aging materials states  {\it the stress depends not only on the instantaneous strain but also on the strain history.} Under this premise, and following  classical Timoshenko's assumptions on a beam filament of length $l>0$, it is rigorously deducted in \cite[Sect. 2]{michele-SIAM}  the following viscoelastic constitutive laws
 \begin{equation}\label{shear-bending-viscoelastic}
	\left\{\begin{array}{l}
\displaystyle		S  \,= \,  k\big(\varphi_{x}+\psi\big) - \int_{0}^{\infty}\mu(s)  \big(\varphi_{x}+\psi\big)(t-s)\,ds , \medskip \\ 
\displaystyle	M\,=\, b\psi_x - \int_{0}^{\infty}g(s) \psi_x(t-s)\,ds,
	\end{array} \right.
\end{equation}
where $\mu,g$ are non-negative relaxation functions, so-called {\it memory kernels}. 
The difference here, when compared to \cite{michele-SIAM}, relies on the fact that we are considering the constitutive laws \eqref{shear-bending-viscoelastic} in the history framework, whereas in \cite[Sect. 2]{michele-SIAM} the analysis is done with null history (see equations (2.7)-(2.8) therein). Some authors also say {\it finite memory} to design memory without history. Hence, with identities \eqref{shear-bending-viscoelastic} in hand, one can obtain at least three different viscoelastic Timoshenko systems, depending on where we consider the viscoelastic couplings \eqref{shear-bending-viscoelastic} in \eqref{tim1}, namely, on the shear force only (cf. \cite{michele-SIAM} with respect to null history) and then a partially damped system emerges; or just on the bending moment (see   \cite{Ammar.Ben.Jaime.Racke.2003,jaime-hug0-jmaa08} for both treatments with or without history) still yielding a partially damped system; or else viscoelastic coupling on both the shear and bending forces by producing a fully damped system (we refer to  \cite{grasseli-pata-prouse-04} for a slightly modified problem with history).

On the other hand, when 
 the beam model is subject to  unknown temperature distribution, then the  principles in thermoelasticity 
 state   {\it the stress depends not only on the elastic strain but also on the thermal strain.} By following up this setting, and still assuming the Timoshenko hypotheses for thin beams, one can find in \cite[Sect. 2]{michele-JEE}  a precise justification of the following thermoelastic constitutive laws
  \begin{equation}\label{shear-bending-therm}
 	\left\{\begin{array}{l}
 		\displaystyle		S  \,= \,  k\big(\varphi_{x} +  \psi\big)  -  \sigma\upsilon  , \medskip \\ 
 		\displaystyle	M\,=\,  b  \psi_{x }  - \varsigma  \vartheta ,
 	\end{array} \right.
 \end{equation}
where $ \sigma, \varsigma>0$ are coefficients related to the thermal expansion, and  $\upsilon=\upsilon(x,t), \,\vartheta=\vartheta(x,t)$   are temperature components  standing for  temperature deviations from a reference state    along 
the longitudinal and vertical directions.  Therefore, in possession of the thermal identities in \eqref{shear-bending-therm}, we can also set at least three distinct thermoelastic Timoshenko systems   depending now where we regard the thermal couplings \eqref{shear-bending-therm} in \eqref{tim1}. 
For instance, with thermal coupling  only on the bending moment, a partially damped system arises (cf. \cite{Rivera.Racke.2002}  where  Fourier's law is taken into account); under the thermal component solely on the shear force we still obtain a partially damped system (as proposed in \cite{jaime-dilberto-santos-2014-1} again under Fourier's law); and last by invoking the thermal coupling on both the bending moment and the shear force we arrive at a fully damped thermoelastic Timoshenko system (as studied in \cite{michele-JEE}). All these possibilities are illustrated in \cite[Sect. 2]{michele-JEE}, as well as the existing literature
dealing with other thermal laws for the heat flux of conduction  is provided  therein (e.g. Gurtin-Pipkin, 
Maxwell-Cattaneo, 
Coleman-Gurtin, 
type III).

Under the above statements, one sees that   the stability of \eqref{tim1}  ultimately depends upon damping's feedback provided by the viscoelastic
coupling \eqref{shear-bending-viscoelastic} or the thermoelastic coupling \eqref{shear-bending-therm} on both (or not) forces of the system. We refer, for instance, the  stability results in
\cite{jaime-dilberto-santos-2014-1,michele-JEE,michele-SIAM,alves-marcio-matofu-rivera-zamp,alves-marcio-matofu-rivera,Ammar.Ben.Jaime.Racke.2003,conti,delloro-pata-jde,grasseli-pata-prouse-04,jaime-hug0-jmaa08,Rivera.Racke.2002,dilberto-mauro,santos-jde2012} just to name few. Additionally, by laying down possible  {\it hybrid} dissipative models generated by 
\eqref{shear-bending-viscoelastic} and \eqref{shear-bending-therm} simultaneously, that is, models featured by mixed damping feedback  in thermo-(visco-)elasticity, then new perspectives  pop up
 in what concerns the stability of Timoshenko systems with temperature and memory terms, as one can see, for example, in   \cite{apalara2018,Feng-aa,hugo-racke,jorge-racke,salim-fareh-2011,salim-fareh-2013,dilberto-mauro,soufyane2008}. In this way, among all possibilities,  we are going to take into account  the  following two ones thermo-(visco-)elastic laws:
\begin{equation} \label{mixed}
\mbox{\texttt{I. Bending and Shear Coupling:} } \ 	\left \{
	\begin{array}{ll}
		S= k\big(\varphi_x +\psi\big) -\sigma\upsilon,  \medskip \\  
		\displaystyle
		M= b\psi_x - \int_{0}^{\infty}{g(s) \psi_x(t-s)\,ds}, 
	\end{array}
	\right.
\end{equation}
or
\begin{equation} \label{mixed-a}
	\mbox{\texttt{II. Bending  Coupling Only:} } \  	\left \{
	\begin{array}{ll}
		S= k\big(\varphi_x +\psi\big),  \medskip \\  
		\displaystyle
		M= b\psi_x - \int_{0}^{\infty}{g(s) \psi_x(t-s)\,ds}-\sigma\upsilon. 
	\end{array}
	\right.
\end{equation}

In what follows, let us first work with  \eqref{mixed}, and then we check what happens to  \eqref{mixed-a}.

\smallskip 

\noindent \texttt{Case I.}
Replacing \eqref{mixed} in \eqref{tim1} we arrive at 
\begin{equation} \label{tim2}
	\left \{
	\begin{array}{l}
	\displaystyle	\rho_1 \varphi_{tt} - k(\varphi_x +\psi)_x + \sigma \upsilon_x= 0  \   \mbox{ in } \  (0,l)\times (0,\infty) ,\medskip \\  
		\displaystyle
		\rho_2\psi_{tt} - b\psi_{xx} + \int_{0}^{\infty}{g(s)\psi_{xx}(t-s)\,ds} + k(\varphi_x + \psi) - \sigma \upsilon= 0  \  \mbox{ in } \  (0,l)\times (0,\infty).
	\end{array}
	\right.
\end{equation} 
Concerning the  temperature deviation, we must provide an equation for the heat flux of conduction. Since  the coupling for this variable is given on the shear force  \eqref{mixed}$_1$, then relying on the facts presented in \cite[Sect. 2]{michele-JEE} we can  derived the following motion equation
\begin{equation} \label{heat1}
\rho_3 \upsilon_t= -q_x - \sigma(\varphi_x + \psi)_t \quad  \mbox{in} \quad  (0,l)\times (0,\infty),
 \end{equation}
where $\rho_3>0$ is a constant related to the heat capacity and $q=q(x,t)$ stands for the heat flux. The Fourier and  Maxwell-Cattaneo laws for the heat flux   are considered in the recent paper \cite{jorge-racke}. Here,  we follow the  Green and Naghdi theory, cf. \cite{gree1,gree2}, to consider the so-called {\it type III thermoelasticity} for the heat flux of conduction, namely,  
\begin{equation} \label{GN}
	q=-\delta p_x - \gamma p_{xt} \quad \mbox{with} \quad p_t=\upsilon,
\end{equation}
where $\delta, \gamma>0$ are constants related to the thermal conductivity, and $p$  stands for  the thermal displacement whose time derivative is  empirically the temperature as in \eqref{GN}. 
Combining  \eqref{heat1} and \eqref{GN}  and plugging the resulting expression into \eqref{tim2}, we obtain the following type III thermoelastic Timoshenko system with long  memory (history)
 \begin{equation} \label{tim3}
 	\left \{
 	\begin{array}{l}
 		\displaystyle	\rho_1 \varphi_{tt} - k(\varphi_x +\psi)_x + \sigma \upsilon_x= 0  \quad  \mbox{in} \quad  (0,l)\times (0,\infty) ,\medskip \\  
 		\displaystyle
 		\rho_2\psi_{tt} - b\psi_{xx} + \int_{0}^{\infty}{g(s)\psi_{xx}(t-s)\,ds} + k(\varphi_x + \psi) - \sigma \upsilon= 0  \quad  \mbox{in} \quad  (0,l)\times (0,\infty),\medskip \\ 
 		\rho_3 \upsilon_{tt}  - \delta\upsilon_{xx}  - \gamma \upsilon_{xxt} + \sigma (\varphi_{x}+\psi)_{tt} \, = \, 0	\quad  \mbox{in} \quad  (0,l)\times (0,\infty).
 	\end{array}
 	\right.
 \end{equation}

For the sake of completeness, we consider  \eqref{tim3} with initial  conditions
\begin{equation} \label{init-cond-v}
	\left \{
	\begin{array}{l}
		\varphi(x,0)= \varphi_0(x), \ \varphi_t(x,0)= \varphi_1(x), \smallskip \\
		\psi(x,s)= \psi_0(x,s), \ \psi_t(x,0)= \partial_t \psi(x,t)|_{t=0}:= \psi_1(x), \ s\leq 0, \smallskip \\
		\upsilon(x,0)= \upsilon_0(x), \upsilon_t(x,0)= \upsilon_1(x), \ x\in(0,l),
	\end{array}
	\right.
\end{equation}
and either boundary conditions 
\begin{subequations}\label{eq271}
	\begin{equation} \label{eq27}
		\left \{
		\begin{array}{lcl}
			\varphi(x,t)= \upsilon(x,t)= 0  & \mbox{for} & x= 0, l; \ t\geq0, \medskip\\
			\psi(x,t)=0   & \mbox{for} &  x=0,l; \ t\in \mathbb{R},
		\end{array}
		\right.
	\end{equation}
	or
	\begin{equation} \label{eq28}
		\left \{
		\begin{array}{lcl}
			\varphi_x(x,t)= \upsilon(x,t)= 0  & \mbox{for} & x= 0, l; \ t\geq0, \medskip\\
			\psi(x,t)=0  & \mbox{for} & x=0,l; \ t\in \mathbb{R}.
		\end{array}
		\right.
	\end{equation}
\end{subequations}

\begin{rem}\label{rem-strength}\rm 
As one can see from the above construction up to achieve \eqref{tim3},  the thermoelastic damping feedback comes from the coupling on the shear force, whereas the viscoelastic dissipative mechanism  is hidden by the memory component coupled on the bending moment, as conducted by \eqref{mixed}. It means that we will be able to extract the strength of both dissipations when dealing with the stability of the solution to \eqref{tim3}. 
In conclusion, we have a fully damped system  not only from the physical point of view but also from the mathematical one, as will be also proved in Section \ref{sec-thermovisco-problem}. 
\end{rem}

In what follows, we introduce two new variables in order to see  problem \eqref{tim3}-\eqref{eq271} in a dissipative and autonomous scenario, namely, as given in \eqref{intro 1.5}-\eqref{bc}.

\noindent {\it Auxiliary temperature variable.} By following a similar idea as introduced in \cite[Sect. 1]{zhang-zuazua} (see also \cite{fatori-jaime-monteiro,dilberto-mauro}) we set the new variable concerning the temperature distribution
\begin{equation}\label{eq29}
	\theta(x,t):= \int_{0}^{t}{\upsilon(x,s)ds} + \dfrac{1}{\delta} z(x),
\end{equation}
where $z\in H_0^1(0,l)$ is the solution of the Cauchy problem
\begin{equation} \label{eq31}
	\left \{
	\begin{array}{ll}
		z_{xx}= \rho_3 \upsilon_1 - \gamma \upsilon_{0,xx} + \sigma (\varphi_{1,x} + \psi_1), & x\in (0,l), \medskip\\
		z(x)=0, & x= 0, l.
	\end{array}
	\right.
\end{equation}

Note that we can formally  write down \eqref{tim3}$_3$ as 
\begin{equation*}\label{eq31a}
	-\delta\left\{\int_{0}^{t}{\upsilon_{xx}(\cdot, s)ds} + \dfrac{1}{\delta} [\rho_3 \upsilon_1 - \gamma \upsilon_{0,xx} + \sigma(\varphi_{1,x} + \psi_1)]\right\} + \rho_3 \upsilon_t - \gamma \upsilon_{xx} + \sigma (\varphi_x + \psi)_t=0,
\end{equation*}
and from \eqref{eq29}-\eqref{eq31} the latter turns into
\begin{equation}\label{eq32}
	\rho_3 \theta_{tt} - \delta \theta_{xx} - \gamma \theta_{xxt} + \sigma (\varphi_x + \psi)_t=0.
\end{equation}

 \noindent {\it  Relative displacement history.} 
 Now, as  introduced by Dafermos   \cite{dafermos} (see also \cite[Sect. 2]{grasseli-pata02})  
 we consider the relative displacement history with respect to the angle rotation
\begin{equation}\label{relativ-hist}
	\eta= \eta (x,t,s):= \psi(x,t)- \psi(x,t-s), \ \  x\in (0,l),\, t\geq 0, \,s> 0.
\end{equation} 
Thus, through \eqref{relativ-hist} the equation \eqref{tim3}$_2$ can be rewritten as 
 \begin{equation}\label{eq300}
  	\rho_2 \psi_{tt} - \left(b - \int_{0}^{\infty}{g(s)ds}\right) \psi_{xx} - \int_{0}^{\infty}{g(s)\eta_{xx}(s)ds} + k(\varphi_x+ \psi) - \sigma \upsilon= 0.
 \end{equation} 
From \eqref{eq300} one sees that a supplementary equation with respect to $\eta$ is necessary to draw up the whole problem. To this purpose, we use again \eqref{relativ-hist} and initial-boundary conditions \eqref{init-cond-v}-\eqref{eq271}
with respect to $\psi$ to  derive formally the next identities
\begin{equation}\label{eq33} 
	\left \{
	\begin{array}{ll}
		\eta_t+ \eta_s= \psi_t & \mbox{in } \ (0,l)\times (0,\infty)\times (0,\infty), \smallskip \\
		\eta(x,t,s)=0 &   \mbox{for } \ x= 0,l; \, t\geq 0, \, s>0, \smallskip \\
		\eta(x,t,0):= \displaystyle \lim_{s\rightarrow 0} \eta(x,t,s)= 0 & \mbox{for } \ x\in (0,l), \, t\geq 0, \smallskip \\
		\eta(x,0,s)= \psi_0(x) - \psi_0(x, -s):= \eta_0(x,s) & \mbox{for } \ x\in (0,l), \, s> 0.
	\end{array}
	\right.
\end{equation}

Therefore, using 
\eqref{eq29}-\eqref{eq33}, and  denoting  $\beta:=  b - \int_{0}^{\infty}{g(s)ds}> 0,$ we can rewrite \eqref{tim3}-\eqref{eq271}
as the following equivalent IBVP
\begin{equation} \label{P1}
	\left \{
	\begin{array}{ll}
		\rho_1\varphi_{tt} - k(\varphi_x + \psi)_x + \sigma \theta_{tx}= 0 & \mbox{in } \ (0,l)\times (0,\infty), \smallskip \\
		\displaystyle\rho_2 \psi_{tt} - \beta \psi_{xx} - \int_{0}^{\infty}{g(s)\eta_{xx}(s)ds} + k(\varphi_x+ \psi) - \sigma \theta_t= 0 & \mbox{in } \ (0,l)\times (0,\infty), \smallskip \\
		\rho_3 \theta_{tt} - \delta \theta_{xx} - \gamma \theta_{xxt} + \sigma (\varphi_x + \psi)_t= 0 & \mbox{in } \ (0,l)\times (0,\infty), \smallskip \\
		\eta_t + \eta_s- \psi_t=0 & \mbox{in } \ (0,l)\times (0,\infty)\times (0,\infty), \smallskip \\
	(\varphi,\varphi_t,\psi,\psi_t,	\theta,\theta_t)(\cdot,0)=(\varphi_0, \varphi_1,\psi_0, \psi_1, \theta_0, \theta_1)(\cdot)  & 	\mbox{in } \ (0,l),\smallskip \\
\eta(\cdot,0,s)= \eta_0(\cdot,s), \ 	\eta(\cdot,t,0)=0  & \mbox{in } \ (0,l), \, s>0,
	\, t\geq 0,	\smallskip \\
\varphi= \psi= \theta=\eta(\cdot,\cdot,s)=0 \ \mbox{or} \ \varphi_x= \psi= \theta=\eta(\cdot,\cdot,s)=0 &  \mbox{on } \, \{0,l\}\times[0,\infty), \, s>0,
	\end{array}
	\right.
\end{equation}
where the relationship between 
 $(\theta_0, \theta_1)$  and $(\upsilon_0, \upsilon_1)$ is assumed as follows
\begin{equation*} 	\theta_1= \upsilon_0 \ \mbox{ and } \
		\rho_3\upsilon_1= \delta \theta_{0,xx} + \gamma \theta_{1,xx} - \sigma (\varphi_{1,x} + \psi_1).
\end{equation*}

Finally, we observe that \eqref{P1} corresponds precisely to problem \eqref{intro 1.5}-\eqref{bc}, which is the main object of study in this article.

\smallskip 

\noindent \texttt{Case II.}   In this case, we replace \eqref{mixed-a} in \eqref{tim1} to obtain
\begin{equation} \label{tim2-a}
	\left \{
	\begin{array}{l}
		\displaystyle	\rho_1 \varphi_{tt} - k(\varphi_x +\psi)_x = 0  \quad  \mbox{in} \quad  (0,l)\times (0,\infty) ,\medskip \\  
		\displaystyle
		\rho_2\psi_{tt} - b\psi_{xx} + \int_{0}^{\infty}{g(s)\psi_{xx}(t-s)\,ds} + k(\varphi_x + \psi) + \sigma \upsilon_x= 0  \quad  \mbox{in} \quad  (0,l)\times (0,\infty).
	\end{array}
	\right.
\end{equation} 
With respect to
the heat flux equation, instead of \eqref{heat1}  we have now (cf. \cite[Sect. 2]{michele-JEE})
\begin{equation} \label{heat1-a}
	\rho_3 \upsilon_t= -q_x - \sigma\psi_{xt} \quad  \mbox{in} \quad  (0,l)\times (0,\infty).
\end{equation}
 
Combining  \eqref{heat1-a} with \eqref{GN}  and substituting 
the resulting expression in  \eqref{tim2-a}, we  obtain the other type III thermoelastic Timoshenko system with history 
\begin{equation} \label{tim3-a}
	\left \{
	\begin{array}{l}
		\displaystyle	\rho_1 \varphi_{tt} - k(\varphi_x +\psi)_x = 0  \quad  \mbox{in} \quad  (0,l)\times (0,\infty) ,\medskip \\  
		\displaystyle
		\rho_2\psi_{tt} - b\psi_{xx} + \int_{0}^{\infty}{g(s)\psi_{xx}(t-s)\,ds} + k(\varphi_x + \psi) + \sigma \upsilon_x= 0  \quad  \mbox{in} \quad  (0,l)\times (0,\infty),\medskip \\ 
		\rho_3 \upsilon_{tt}  - \delta\upsilon_{xx}  - \gamma \upsilon_{xxt} + \sigma \psi_{xtt} \, = \, 0	\quad  \mbox{in} \quad  (0,l)\times (0,\infty).
	\end{array}
	\right.
\end{equation}

We note that system \eqref{tim3-a} is precisely the one studied in \cite{salim-belkacem-09}, see problem (1.8) therein.

Last, proceeding similarly (with minor modifications) as in \eqref{eq29}-\eqref{eq33}, one can formally convert \eqref{tim3-a} into the following autonomous equivalent system
 \begin{equation} \label{P1-a}
 	\left \{
 	\begin{array}{ll}
 		\rho_1\varphi_{tt} - k(\varphi_x + \psi)_x = 0 & \mbox{in } \ (0,l)\times (0,\infty), \smallskip \\
 		\displaystyle\rho_2 \psi_{tt} - \beta \psi_{xx} - \int_{0}^{\infty}{g(s)\eta_{xx}(s)ds} + k(\varphi_x+ \psi) + \sigma \theta_{tx}= 0 & \mbox{in } \ (0,l)\times (0,\infty), \smallskip \\
 		\rho_3 \theta_{tt} - \delta \theta_{xx} - \gamma \theta_{xxt} + \sigma  \psi_{xt}= 0 & \mbox{in } \ (0,l)\times (0,\infty), \smallskip \\
 		\eta_t + \eta_s- \psi_t=0 & \mbox{in } \ (0,l)\times (0,\infty)\times (0,\infty),
 	\end{array}
 	\right.
 \end{equation}
with proper initial-boundary conditions, which corresponds exactly to problem  \eqref{intro 1.5-n}-\eqref{intro 1.8-n}.

 \begin{rem}\label{rem-strength-a}\rm 
It is worth mentioning  that the above construction to reach \eqref{tim3-a}, and consequently \eqref{P1-a},  reveals us   that the thermal and viscoelastic couplings are both given on the bending moment, here conducted by  \eqref{mixed-a}. It means that we can no longer expect  a fully dissipative mechanism. On the contrary, both thermal and viscoelastic damping terms now propagate to the bending moment, and only to the shear force by means of the equal wave speeds assumption \eqref{ews}, as already aforementioned in the results of \cite{salim-belkacem-09}. In conclusion, \eqref{P1-a} represents a partially damped system which  is very different,  from the stability point of view, when compared to  the fully damped problem \eqref{P1},  although {\it close} to it in terms of variables and other stuff. The coming mathematical results certify these formal statements that have been physically built.
 \end{rem}

\section{Main result on stability}\label{sec-thermovisco-problem}	
 
In this section, our main goal is to prove  the exponential stability of  
solutions to problem \eqref{intro 1.5}-\eqref{bc}.  Before doing so, let us first introduce the semigroup setting that will be useful hereafter.

\subsection{Semigroup solution}\label{subsec-existence-thermovisco}

For each boundary condition in \eqref{bc}, we need to take different phase spaces. Here, we consider
\begin{equation*}
  \mathcal{H}_1=
H^1_0(0,l) \times L^2(0,l) \times H_0^1(0,l) \times L^2(0,l)\times H_0^1(0,l)\times L^2(0,l)\times \mathcal{M} \ \ \mbox{ for } \ \ (\ref{bc-a}), 
\end{equation*}
and
\begin{equation*}
\mathcal{H}_2= H^1_*(0,l) \times L^2_*(0,l) \times H_0^1(0,l) \times L^2(0,l) \times H_0^1(0,l) \times L^2(0,l)\times \mathcal{M} \ \ \mbox{ for } \ \ (\ref{bc-b}),
\end{equation*}
where 
$
H^1_*(0,l) = H^1(0,l)\cap L^2_*(0,l), \; L^2_*(0,l)=\left\{u\in L^2(0,l); \ \frac{1}{l}\int_0^l u(x)\,dx=0  \right\}
$ and
$$ \mathcal{M}= \left \{ \eta: (0,\infty) \rightarrow H_0^1(0,l); \ \ \int_0^{\infty} g(s) \|\eta(s)\|_{H_0^1(0,l)}^2 ds < \infty \right \}.$$

It is well-known that  ${{\cal H}_j},$ for each $ j=1,2,$ is a Hilbert space endowed with norm
\begin{align*}%\label{norm-H}
\|U\|^2_{\mathcal{H}_j}
\, = \, \int_0^l \Big[\rho_1|\Phi|^2 +\rho_2|
\Psi|^2 + \rho_3|\Theta_x|^2 +\beta|\psi_x|^2+
k |\varphi_x + \psi|^2+\delta|\theta_x|^2  + \int_{0}^{\infty}{g(s) |\eta_x(s)|^2 ds}\Big]dx,
\end{align*}
for $U=( \varphi, \Phi, \psi, \Psi,\theta, \Theta, \eta) \in \mathcal{H}_j,$ and respective scalar product   $(\cdot,\cdot)_{\mathcal{H}_j}$.

\smallskip

Under the above notation and setting $\Phi:=\varphi_t, \Psi:=\psi_t, \Theta:=\theta_t$,  we can
 convert the particular system (\ref{intro 1.5})-(\ref{bc}) into the following abstract problem
\begin{equation} \label{semigroup 3.1}
\left\{\begin{array}{l}
U_t = {\cal A}_j\, U, \quad  t>0, \medskip \\
U(0) =(\varphi_0, \varphi_1, \psi_0, \psi_1,\theta_0, \theta_1, \eta_0):=U_0, \end{array}\right.
\end{equation}
where $\mathcal{A}_j:D({\cal A}_j)\subset \mathcal{H}_j\rightarrow \mathcal{H}_j$
is defined by
\begin{equation} \label{semigroup 3.2}
{\cal A}_jU = \left( 
\begin{array}{c}
 \Phi   \\ 
\displaystyle \frac{k}{\rho_1}(\varphi_x+\psi)_x- \frac{\sigma}{\rho_1}\,\Theta_x  \medskip  \\ 
\Psi   \\
\dfrac{1}{\rho_2} \left( \beta \psi + \int_{0}^{\infty}{g(s)\eta(s)ds}\right)_{xx} - \dfrac{k}{\rho_2}(\varphi_x+\psi)   + \dfrac{\sigma}{\rho_2} \Theta   \medskip  \\ 
\Theta  \\
\displaystyle  \frac{1}{\rho_3}(\delta \theta + \gamma \Theta)_{xx}-\frac{\sigma}{\rho_3}(\Phi_x+\Psi)\\
\Psi - \eta_s
\end{array} \right), 
\end{equation}
for any $U=( \varphi, \Phi, \psi, \Psi,\theta, \Theta, \eta) \in D({\cal A}_j),$ with domain
\begin{align*}
D(\mathcal{A}_1)= \bigg\{U\in \mathcal{H}_1 \, | & \, \Phi, \Psi, \Theta\in H_0^1(0,l), \, \eta_s \in \mathcal{M}, \eta(0)=0,\\
& \ \varphi,\delta\theta+\gamma\Theta, \beta \psi+  \int_{0}^{\infty}{g(s)\eta(s)ds} \in H^2(0,l)\bigg\} \qquad \mbox{for} \qquad \quad \eqref{bc-a}
\end{align*}
and
\begin{align*}
    D(\mathcal{A}_2)= \bigg\{U\in \mathcal{H}_2 \, | & \, \Phi\in H_*^1(0,l), \varphi_x, \Psi, \Theta \in H_0^1(0,l), \eta_s \in \mathcal{M}, \eta(0)=0,\\
& \ \varphi,\delta\theta+\gamma\Theta, \beta \psi+  \int_{0}^{\infty}{g(s)\eta(s)ds} \in H^2(0,l)\bigg\}\qquad  \mbox{for} \qquad \quad \eqref{bc-b}.
\end{align*}

Under the above construction, it is easy to verify that $0 \in \rho(\mathcal{A}_j),\, j=1,2$, and also ${\cal A}_j$ is dissipative on ${{\cal H}_j}$ with 
\begin{equation}\label{diss1}
\mbox{Re }({\cal A}U,U)_{{{\cal H}_j}}=
-\gamma\int_0^l|\Theta_x|^2 dx  + \frac{1}{2}\int_{0}^{\infty}{g'(s) \int_{0}^{l}{|\eta_x(s)|^2dx}\, ds}, \ \ \forall \,  U\in D({\cal A}_j), \ j=1,2.
\end{equation}
To perform integration by parts in the last term of \eqref{diss1}, we proceed analogously to \cite{grasseli-pata02}, see pages 162-163 therein.

Therefore, relying on  Pazy's book \cite{pazy}, we have that ${\cal A}_j$ is the infinitesimal generator of a $C_0$-semigroup of contractions $T(t)=e^{\mathcal{A}_jt}$ on $\mathcal{H}_j,$ and existence and uniqueness result reads as follows. 
\begin{Theorem}\label{theo-existence} Under the above notations and Assumption $\ref{Hyp-g}$ we have for $j=1,2:$ 
	\begin{itemize}
		\item[$ (i) $] If $U_0\in{{\cal H}_j},$ then problem  $ (\ref{semigroup 3.1}) $ has a unique mild solution $U\in C^0([0,\infty),{{\cal H}_j})$ given by 
		$$U(t)=e^{\mathcal{A}_jt}U_0, \quad t\geq0.$$
		
		\item[$ (ii) $]  If $U_0\in D({\cal A}_j),$  then problem  $ (\ref{semigroup 3.1}) $ has a unique regular solution
		$$
		U\in C^0([0,\infty),D({\cal A}_j))\cap  C^{1}([0,\infty),{\cal{H}}_j).
		$$
		\item[$ (iii) $]  If $U_0\in D({{\cal A}_j}^n),$  $n\geq2$ integer, then the solution is more regular 
		$$
		U\in
		\bigcap_{\nu=0}^{n}C^{n-\nu}([0,\infty),D({\cal A}^{\nu})).
		$$
	\end{itemize}
\end{Theorem}

\subsection{Exponential stability}\label{subsec-expdecay}

Our main result in this section is
 given below.

\begin{Theorem}[Exponential Stability]\label{theo-expo3} 
	Under the above notations  and Assumption $\ref{Hyp-g}$, there exist constants 	$M,m>0$ independent of $U_0\in \mathcal{H}_j$ such that the semigroup solution $U(t)= e^{\mathcal{A}_jt}U_0$ decays as
	\begin{equation*}\label{exp-decay3}
	\|U(t)\|_{\mathcal{H}_j} \leq  M e^{-m \, t} \|U_0\|_{\mathcal{H}_j}, \quad t>0,\; j=1, 2.
	\end{equation*}
In other words, the thermo-viscoelastic system \eqref{intro 1.5}-\eqref{bc} is exponentially stable. 	
\end{Theorem}
 
The proof of Theorem \ref{theo-expo3} is based on the well-known characterization of exponential stability for $C_0$-semigroups of contractions, cf. \cite{gearhart,huang,pruss}. See also  \cite[Thm. 1.3.2]{Zhuangyi.Liu.book}. Accordingly, we need to proof the following properties:
\begin{equation}\label{properties}
i \mathbb{R}\subset \rho(\mathcal{A}_j) \ \ \mbox{ and } \ \ \limsup_{|\lambda|\to\infty}	\| (i \lambda I_d - \mathcal{A}_j)^{-1}  \|_{{\cal L}(\mathcal{H}_j)}<\infty, \ \ j=1,2,
\end{equation}
where $ \mathcal{A}_j$ is defined in \eqref{semigroup 3.2}.

\subsubsection{Resolvent set}

Let us assume by contradiction that $i \mathbb{R}\not\subset \rho(\mathcal{A}_j), \, j=1,2.$ Appealing to  Proposition \ref{result-principal-apA},
there exist a constant $\omega\in (0, \ell]$, $\ell=\|(-\mathcal{A}_j)^{-1}\|_{\mathcal{L}(\mathcal{H}_j)}^{-1}$, a sequence $\xi_n\in \mathbb{R}$, with $\xi_n \rightarrow w$ and $|\xi_n|< w$, and a sequence of vector functions
\begin{equation}\label{conv-norma}
U_n= (\varphi_n, \Phi_n, \psi_n, \Psi_n, \theta_n, \Theta_n, \eta_n) \in D(\mathcal{A}_j)\quad \mbox{with}\quad \|U_n\|_{\mathcal{H}_j}=1,
\end{equation}
such that
 \begin{equation} \label{e-seq-temp-resolv1234}
i \xi_n U_n - \mathcal{A}_j U_n \rightarrow 0\quad \mbox{in}\quad \mathcal{H}_j,\; j=1, 2.
\end{equation}
Using the expression for  $ \mathcal{A}_j$ given in \eqref{semigroup 3.2}, then \eqref{e-seq-temp-resolv1234} can be rewritten in terms of its components as follows
\begin{equation}\label{elep-regul}
\left\{\begin{array}{lcl}
i\xi_n \varphi_n - \Phi_n\rightarrow 0 & \mbox{ in } & H_{0}^{1}(0,l)\; \mbox{or}\; H_{*}^{1}(0,l), \medskip  \\
i\xi_n\rho_1 \Phi_n - k(\varphi_{n, x}+\psi_n)_{x} + \sigma \Theta_{n,x} \rightarrow 0 & \mbox{ in } & L^2(0,l)\; \; \mbox{or}\; L^2_*(0,l),
\medskip  \\
i \xi_n \psi_n-\Psi_n \rightarrow 0  & \mbox{ in } & H_{0}^{1}(0,l),
\medskip  \\
i\xi_n \rho_2\Psi_n + k(\varphi_{n,x}+ \psi_n)   & & \\
\qquad \quad \; \; - \left( \beta\psi_{n} + \int_{0}^{\infty}{g(s)\eta_{n}(s)ds} \right)_{xx} - \sigma \Theta_n\rightarrow 0 & \mbox{ in } & L^{2}(0,l), \medskip  \\
i \xi_n \theta_n - \Theta_n \rightarrow 0  & \mbox{ in } & H_{0}^{1}(0,l),\medskip  \\
i \xi_n \rho_3 \Theta_n - (\delta \theta_n + \gamma \Theta_n)_{xx} + \sigma (\Phi_{n,x} + \Psi_n) \rightarrow 0 & \mbox{in} & L^2(0,l),\medskip  \\
i\xi_n \eta_n + \eta_{n,s} -\Psi_n \rightarrow 0 &\mbox{in} & \mathcal{M}.
\end{array}\right.
\end{equation}

Now, our purpose is to prove that 
\begin{equation}\label{conv-to-zero}
\|U_n\|^2_{\mathcal{H}_j} \rightarrow 0, \ \ j=1,2,
\end{equation}
which provides the desired contradiction with \eqref{conv-norma} and, therefore,  the proof of 
$i \mathbb{R}\subset \rho(\mathcal{A}_j), \, j=1,2,$ is complete.

 In what follows, the proof of \eqref{conv-to-zero} will be done through some lemmas, where due to the nature of the boundary conditions \eqref{bc-a} and \eqref{bc-b} we first conclude it for $U_n\in   D(\mathcal{A}_2), j=2,$ and then for  $U_n\in   D(\mathcal{A}_1), j=1$. 

\begin{Lemma}\label{e-seq-tem-resultado1-compl}
Under the assumptions of Theorem $\ref{theo-expo3}$ and the above notations, let $U_n$ be the sequence satisfying \eqref{e-seq-temp-resolv1234}. Then, 
	\begin{enumerate}
		\item $\Theta_n, \; \theta_n \rightarrow 0   \ \ \mbox{in} \ \ H_0^1(0,l)$, as $n\rightarrow \infty$; \smallskip
		
		\item $\displaystyle\int_{0}^{\infty}{[- g'(s)] \|\eta_{n, x}(s)\|^2_{L^2}ds}  \rightarrow 0$, as $n\rightarrow \infty$;
		\smallskip
		
		\item $ \eta_n\rightarrow 0  \ \ \mbox{in} \ \  \mathcal{M}$, as $n\rightarrow \infty$.
	\end{enumerate}
\end{Lemma}
\begin{proof}
	From \eqref{diss1} and \eqref{e-seq-temp-resolv1234} we promptly get 
	\begin{equation*}
	\gamma \|\Theta_{n,x}\|_{L^2}^2 + \dfrac{1}{2} \int_{0}^{\infty}{[-g'(s)] \|\eta_{n, x}(s)\|^2_{L^2}ds}= \re(i \xi_n U_n -\mathcal{A}_j U_n, U_n)_{\mathcal{H}_j}\rightarrow 0.
	\end{equation*}
From this, 	\eqref{elep-regul}$_5$ and condition \eqref{sobreg}, the three statements of Lemma \ref{e-seq-tem-resultado1-compl} are fulfilled.
\end{proof}

Under the limits of Lemma \ref{e-seq-tem-resultado1-compl}, the convergences of \eqref{elep-regul} can be reduced into the next ones
\begin{subequations}\label{ee1-seq-compl2}
	\begin{empheq}{align}
	 i\xi_n \varphi_n - \Phi_n\rightarrow 0  & \ \mbox{ in } \ H_{0}^{1}(0,l)\; \mbox{or}\; H_{*}^{1}(0,l),& \label{ee11-seq-compl2} \medskip \\
 i\xi_n\rho_1 \Phi_n - k(\varphi_{n, x}+\psi_n)_{x} \rightarrow 0  &  \ \mbox{ in } \ L^2(0,l)\; \, \mbox{or}\; L^2_*(0,l),& \label{ee12-seq-compl2} \medskip \\
 i \xi_n \psi_n-\Psi_n \rightarrow 0  & \ \mbox{ in } \ H_{0}^{1}(0,l),& \label{ee13-seq-compl2}\medskip \\
	i\xi_n \rho_2\Psi_n + k(\varphi_{n,x}+ \psi_n)  - \left( \beta\psi_{n} + \int_{0}^{\infty}{g(s)\eta_{n}(s)ds} \right)_{xx} \rightarrow 0  & \ \mbox{ in } \ L^2(0,l), & \label{ee14-seq-compl2} \medskip \\
	i\xi_n \eta_n + \eta_{n,s} -\Psi_n \rightarrow 0   &  \ \mbox{ in } \ \mathcal{M}.& \label{ee17-seq-compl2}
	\end{empheq}
\end{subequations}

\begin{Lemma}\label{e-seq-temp-resultado3}
Under the assumptions of Theorem $\ref{theo-expo3}$ and the above notations, let $U_n$ be the sequence satisfying \eqref{e-seq-temp-resolv1234}. Then,
	\begin{equation}\label{e-seq-temp-resultado3-equation}
	\Psi_n, \psi_n \rightarrow 0\quad \mbox{in}\quad H_{0}^{1}(0,l).
	\end{equation}
\end{Lemma}
\begin{proof} 
Firstly, from \eqref{ee13-seq-compl2} we see that
\begin{equation}\label{Psi_nx-conv-zero-1}
   i \xi_n(\psi_{n,x}, \Psi_{n,x})_{L^2} - \|\Psi_{n,x}\|_{L^2}^2 \rightarrow 0.
\end{equation}
Using Cauchy-Schwarz and Young's inequalities, we get
\begin{equation}\label{Psi_nx-conv-zero-3}
    \|\Psi_{n,x}\|_{L^2}^2\leq 2\big|i\xi_n (\psi_{n,x}, \Psi_{n,x})_{L^2} - \|\Psi_{n,x}\|_{L^2}^2\big| + \xi_n^2 \|\psi_{n,x}\|_{L^2}^2.
\end{equation}
Thus, from 
 \eqref{Psi_nx-conv-zero-1}-\eqref{Psi_nx-conv-zero-3}, since  $\xi_n$ is bounded and $\|\psi_{n,x}\|^2_{L^2}\leq \frac{1}{\beta}\|U_n\|^2_{\mathcal{H}_j}$,  we obtain that $(\|\Psi_{n,x}\|_{L^2})_{n\in \mathbb{N}}$ is bounded.

On the other hand, since $\eta_n\in \mathcal{M}$, we have  $g\|\eta_x(\cdot)\|_{L^2}^2\in L^1((0,\infty))$ and also (see again \cite{grasseli-pata02})
    \begin{equation}\label{210}
     \lim_{z\rightarrow \infty} g(z)\|\eta_{n,x}(z)\|_{L^2}^2=0.
    \end{equation}

We can see that $s\mapsto \frac{1}{\xi_n} \Psi_n \in \mathcal{M}$ happens for any $n\in \mathbb{N}$. Now, taking the multiplier $\frac{1}{\xi_n^2}g(s) \overline{\Psi}_n$ in \eqref{ee17-seq-compl2} and integrating on $(0,l)\times(0, \infty)$, we have
\begin{equation}\label{e-seq-temp-resultado37}
    i\left(\eta_n(\cdot), \frac{\Psi_n}{\xi^2_n} \right)_{\mathcal{M}} + \underbrace{\frac{1}{\xi_n^2} \left(\eta_{n,s}(\cdot), \Psi_n
\right)_{\mathcal{M}}}_{:=S_n} - \frac{1}{\xi_n^2} \left(\Psi_n, \Psi_n
\right)_{\mathcal{M}} \rightarrow 0.
\end{equation}
We claim that $S_{n}\rightarrow 0$. Indeed, let us fix $n\in \mathbb{N}$. Integrating $S_n$ by parts with respect to $s$, using Lemma \ref{e-seq-tem-resultado1-compl}, the fact that $(\Psi_n)_{n\in \mathbb{N}}$ is bounded in $H_{0}^{1}(0,l)$, \eqref{210}, and Fubini's Theorem, we get
\begin{eqnarray*}
   \nonumber |S_n| &=& \left|-\dfrac{1}{\xi^2_n}\int_{0}^{\infty}{g(s)(\eta_{n,s}(s), \Psi_n)_{H_0^1} ds}\right|\\
   \nonumber &=&
   \left|-\dfrac{1}{\xi^2_n}\int_{0}^{l}{\overline{\Psi_{n,x}}\left(g(s) \eta_{n,x}(s)\bigg|_{0}^{\infty} - \int_{0}^{\infty}{g'(s) \eta_{n,x}(s)ds} \right) dx} \right|\\
  \nonumber &\leq& \dfrac{1}{\xi^2_n} \|\Psi_{n,x}\|_{L^2} \left\|\int_{0}^{\infty}{g'(s)\eta_{n,x}(s) ds} \right\|_{L^2}\\
   &\leq& \dfrac{1}{\xi_n^2}\left(\int_{0}^{\infty}{[-g'(s)]ds} \right)^{\frac{1}{2}} \left( \int_{0}^{\infty}{[-g'(s)] \|\eta_{n,x}(s)\|_{L^2}^2 ds} \right)^{\frac{1}{2}} \|\Psi_{n,x}\|_{L^2} \rightarrow 0.
\end{eqnarray*}
Thus, since $g(0)<\infty$,  \eqref{e-seq-temp-resultado37} and Lemma \ref{e-seq-tem-resultado1-compl} imply  $\Psi_n\rightarrow 0$ in $H_0^1(0,l)$ and, consequently, \eqref{ee13-seq-compl2} yields $\psi_n \rightarrow 0$ in $H_0^1(0,l)$, which proves \eqref{e-seq-temp-resultado3-equation} as desired.
\end{proof}

The next result still holds for both boundary conditions at the same time.

\begin{Lemma}\label{e-seq-temp-resultado5}
Under the assumptions of Theorem $\ref{theo-expo3}$ and the above notations, let $U_n$ be the sequence satisfying \eqref{e-seq-temp-resolv1234}. Then,
	\begin{equation}\label{e-seq-temp-resultado5-eq1}
	-\rho_1 \|\Phi_n\|_{L^2}^2  + k \|\varphi_{n,x} + \psi_n\|_{L^2}^2 \rightarrow 0,
	\end{equation}
	and
	\begin{equation}\label{e-seq-temp-resultado-front}
	-  \left(\beta \psi_{n,x} + \int_{0}^{\infty}{g(s) \eta_{n, x}(s) ds}\right)\cdot \overline{\varphi_{n,x}}\bigg|_{0}^{l}  + k \|\varphi_{n,x} + \psi_n\|_{L^2}^2 \rightarrow 0.
	\end{equation}
\end{Lemma}

\begin{proof}
On the one hand, taking the multipliers $\rho_1 \overline{\Phi_n}$ in \eqref{ee11-seq-compl2} and  $\overline{\varphi_n}$ in \eqref{ee12-seq-compl2} and adding the resulting expressions, we obtain
\begin{equation*}%\label{resol-seq-eresultado31}
    i\rho_1\xi_n [(\varphi_n, \Phi_n)_{L^2} + (\Phi_n, \varphi_n)_{L^2}] - \rho_1 \|\Phi_n\|_{L^2}^2 - k ((\varphi_{n,x} + \psi_n)_x, \varphi_n)_{L^2} \rightarrow 0.
\end{equation*}
Integrating by parts and taking the real part, we have
\begin{equation}\label{resol-seq-eresultado32}
    - \rho_1 \|\Phi_n\|_{L^2}^2 + k \mbox{Re}\, (\varphi_{n,x} + \psi_n, \varphi_{n,x})_{L^2} \rightarrow 0.
\end{equation}

On the other hand, taking the multiplier $\rho_2 \Psi_n$ in \eqref{ee13-seq-compl2} and $\overline{\psi_n}$ in \eqref{ee14-seq-compl2}, and adding the resulting convergences, we get
\begin{eqnarray*}%\label{resol-eq-333}
  \nonumber   i\rho_2 \xi_n [(\psi_n, \Psi_n)_{L^2} + (\Psi_n, \psi_n)_{L^2}] &-& \left(\left( \beta \psi_{n} + \int_{0}^{\infty}{g(s) \eta_n(s)ds} \right)_{xx}, \psi_n \right)_{L^2}\\
    &+&  k (\varphi_{n,x} + \psi_n, \psi_n)_{L^2} - \rho_2 \|\Psi_n\|_{L^2}^2 \rightarrow 0.
\end{eqnarray*}
Performing integration by parts and using Lemmas \ref{e-seq-tem-resultado1-compl} and \ref{e-seq-temp-resultado3},   we infer
\begin{equation}\label{resolv-seq-eresultado33}
  k \mbox{Re}\, (\varphi_{n,x} + \psi_n, \psi_n)_{L^2} \rightarrow 0.
\end{equation}
Computing \eqref{resol-seq-eresultado32}+\eqref{resolv-seq-eresultado33}, we arrive at
\eqref{e-seq-temp-resultado5-eq1}.

%\begin{equation*}\label{resolv-seq-eresultado34}
%  -\rho_1 \|\Phi_n\|_{L^2}^2  + k \|\varphi_{n,x} + \psi_n\|_{L^2}^2 \rightarrow 0,
%\end{equation*}
%which proves 

\smallskip

Now, taking the multiplier $k(\overline{\varphi}_{n,x} + \overline{\psi}_n)$ in \eqref{ee14-seq-compl2}, we can write
\begin{eqnarray*}
  \nonumber  && i\beta_n \rho_2 k (\Psi_n, \varphi_{n,x} + \psi_n)_{L^2} - \left(\left(\beta \psi_{n} + \int_{0}^{\infty}{g(s) \eta_n(s)ds}\right)_{xx}, k (\varphi_{n,x} + \psi_n) \right)_{L^2}\\
    && \qquad \qquad \qquad \qquad \qquad \qquad \qquad \qquad \qquad \qquad  + k^2\|\varphi_{n,x} + \psi_n\|_{L^2}^2 \rightarrow 0.
\end{eqnarray*}
Performing integration by parts, using the limits $\eta_n \rightarrow 0$ in $\mathcal{M}$, $\psi_n, \Psi_n \rightarrow 0$ in $H_0^1(0,l)$ and the boundedness of $(\varphi_{n,x} + \psi_n)_{n\in \mathbb{N}}$  in $L^2(0,l)$, we get
\begin{eqnarray}\label{e-seq-temp-resultado35}
   \nonumber &&- k \left(\beta \psi_{n,x} + \int_{0}^{\infty}{g(s) \eta_{n, x}(s) ds}\right)\cdot (\overline{\varphi_{n,x}} + \overline{\psi_n})\bigg|_{0}^{l}  + k^2 \|\varphi_{n,x} + \psi_n\|_{L^2}^2\\
    &&  + \left(\beta \psi_{n,x} + \int_{0}^{\infty}{g(s) \eta_{n,x}(s) ds}, k (\varphi_{n,x} + \psi_n)_x \right)_{L^2} \rightarrow 0.
\end{eqnarray}
On the other hand, taking the multiplier $\beta \psi_{n,x} + \int_{0}^{\infty}{g(s) \eta_{n,x}(s) ds}$ in \eqref{ee12-seq-compl2}, we see
\begin{eqnarray*}\label{e-seq-temp-resultado38}
   &&-i \rho_1 \xi_n \left(\beta \psi_{n,x}+ \int_{0}^{\infty}{g(s) \eta_{n,x}(s) ds}, \Phi_n \right)_{L^2}\\ 
	&&- \left(\beta \psi_{n,x} + \int_{0}^{\infty}{g(s)  \eta_{n,x}(s) ds}, k (\varphi_{n,x} + \psi_n)_x\right)_{L^2} \rightarrow 0.
\end{eqnarray*}
Using that $(\|\Phi_n\|_{L^2})_{n\in \mathbb{N}}$ is bounded and $\left\|\beta \psi_{n,x} + \int_{0}^{\infty}{g(s) \eta_{n,x}(s)ds} \right\|_{L^2} \rightarrow 0$, we have from the previous limit that
\begin{equation}\label{e-seq-temp-resultado39}
    - \left(\beta \psi_{n,x} + \int_{0}^{\infty}{g(s)  \eta_{n,x}(s) ds}, k (\varphi_{n,x} + \psi_n)_x\right)_{L^2}
\rightarrow 0.
\end{equation}
 Computing \eqref{e-seq-temp-resultado35} + \eqref{e-seq-temp-resultado39}, we obtain \eqref{e-seq-temp-resultado-front}, which completes the proof of Lemma \ref{e-seq-temp-resultado5}.
 \end{proof}

We are finally in the position to conclude the proof of \eqref{conv-to-zero}. We proceed in two cases as follows. 

\smallskip 

\noindent \underline{\it Proof of \eqref{conv-to-zero} for  $j=2.$} 
In this case, we are dealing with boundary condition \eqref{bc-b}. Thus, from  \eqref{e-seq-temp-resultado-front} we directly have $\|\varphi_{n,x} + \psi_n\|_{L^2}\rightarrow 0$, which in turn combined with \eqref{e-seq-temp-resultado5-eq1} leads to $\rho_1 \|\Phi_n\|_{L^2}^2 \rightarrow 0$. From these limits and  Lemmas \ref{e-seq-tem-resultado1-compl}-\ref{e-seq-temp-resultado5}, we conclude  
\begin{equation*}
\|U_n\|^2_{\mathcal{H}_2} \rightarrow 0.
\end{equation*}
\qed

\noindent \underline{\it Proof of \eqref{conv-to-zero} for  $j=1.$} In this case, we are under the boundary condition \eqref{bc-a} and more computations are necessary. Initially, we observe that   Lemmas \ref{e-seq-tem-resultado1-compl}-\ref{e-seq-temp-resultado5} imply that  \eqref{ee11-seq-compl2}-\eqref{ee17-seq-compl2} can be rewritten as
\begin{subequations}\label{ef1}
	\begin{empheq}[left=\empheqlbrace]{align}
	& i\xi_n \varphi_n - \Phi_n\rightarrow 0 \qquad \quad \qquad \qquad \qquad \qquad \quad \;  \, \, \mbox{in}\; \quad H_{0}^{1}(0,l)\; \mbox{or}\;  H_{*}^{1}(0,l), \label{ef11} \medskip\\
	& i\xi_n\rho_1 \Phi_n - k\varphi_{n, xx} \rightarrow 0 \quad \qquad \qquad \qquad \qquad \quad \mbox{in}\; \quad L^2(0,l)\; \, \mbox{or}\; L^2_*(0,l), \label{ef12} \medskip\\
	& k\varphi_{n,x}  - \left( \beta\psi_{n} + \int_{0}^{\infty}{g(s)\eta_{n}(s)ds} \right)_{xx} \rightarrow 0 \quad \;  \; \mbox{in}\; \quad L^2(0,l), \label{ef13} \medskip \\
	&\eta_{n,s}(s) \rightarrow 0\qquad \qquad \qquad \qquad \qquad \qquad \qquad \; \, \, \,  \mbox{in}\; \quad \mathcal{M}. \label{ef15}
	\end{empheq}
\end{subequations} 
 
 Taking the multiplier $\left( x- \frac{l}{2}\right) \left(\beta \overline{\psi_{n,x}}+ \int_{0}^{\infty}{g(s)\overline{\eta_{n,x}} ds}\right)$ in \eqref{ef13}, we have
 \begin{align}\label{lemma-e-seq100}
 \nonumber  &\left(k\varphi_{n,x} \left( x- \frac{l}{2}\right) , \beta \psi_{n,x} + \int_{0}^{\infty}{g(s) \eta_{n,x}(s) ds}  \right)_{L^2}\\ 
 & -\left(\left[\beta\psi_n + \int_{0}^{\infty} g(s) \eta_n(s) ds\right]_{xx} \left( x- \frac{l}{2}\right), \beta \psi_{n,x} + \int_{0}^{\infty}{g(s) \eta_{n,x}(s) ds} \right)_{L^2}\rightarrow 0.
 \end{align}

Using, in \eqref{lemma-e-seq100}, that $\|\beta \psi_{n, x} + \int_{0}^{\infty}{g(s) \eta_{n,x}(s) ds}\|_{L^2}\rightarrow 0$ and $(\|\varphi_{n,x}\|_{L^2})_{n\in \mathbb{N}}$ is bounded,  and taking the real part in resulting convergence, we obtain
\begin{equation*}
    - \re \left(\left[\beta\psi_n + \int_{0}^{\infty} g(s) \eta_n(s) ds\right]_{xx}  \left( x- \frac{l}{2}\right), \beta \psi_{n,x} + \int_{0}^{\infty}{g(s) \eta_{n,x}(s) ds} \right)_{L^2} \rightarrow 0.
\end{equation*}
Now, using $\frac{d}{dx}|u|^2= 2 \re (\overline{u}\, u_x)$ and performing integration by parts, we get
\begin{align*}
 & \dfrac{l}{2} \left(\left|\beta \psi_{n,x}(0) + \int_{0}^{\infty}{g(s) \eta_{n,x}(0,s) ds} \right|^2  + \left|\beta \psi_{n,x}(l) + \int_{0}^{\infty}{g(s) \eta_{n,x}(l,s) ds} \right|^2 \right)\\
  & \qquad \qquad \qquad \qquad \qquad \qquad\qquad \qquad - \bigg \|\beta\psi_{n,x} + \int_{0}^{\infty}{g(s)
 \eta_{n,x}(s) ds} \bigg \|_{L^2}^2 \rightarrow 0.
\end{align*}
Again, using $\|\beta \psi_{n, x} + \int_{0}^{\infty}{g(s) \eta_{n,x}(s) ds}\|_{L^2}\rightarrow 0$, we arrive at
\begin{equation}\label{front-psi+int}
    \beta \psi_{n,x}(0) + \int_{0}^{\infty}{g(s) \eta_{n,x}(0,s) ds}, \quad \beta \psi_{n,x}(l) + \int_{0}^{\infty}{g(s) \eta_{n,x}(l,s) ds} \rightarrow 0.
\end{equation}

On the other hand, taking the multiplier $\left(x-\frac{l}{2}\right) \left(\varphi_{n,x} + \psi_n \right)$ in \eqref{ef12}, we infer \begin{equation}\label{cond-a}
   \rho_1 i \xi_n \left(\Phi_n   \left(x - \frac{l}{2}\right), \varphi_{n,x} + \psi_{n} \right)_{L^2} - k\left((\varphi_{n,x} + \psi_n)_x  \left(x - \frac{l}{2}\right), \varphi_{n,x} + \psi_n \right)_{L^2} \rightarrow 0.
\end{equation}
From convergences \eqref{ee11-seq-compl2} and \eqref{ee13-seq-compl2}, we deduce
\begin{equation}\label{cond-a1}
    -i\xi_n \overline{\left(\varphi_{n,x} + \psi_n\right)} - \overline{\left(\Phi_{n,x} + \Psi_n \right)}\rightarrow 0 \quad \mbox{in}\quad L^2(0,l).
\end{equation}
Taking the multiplier $\rho_1 \Phi_n  \left(x - \frac{l}{2}\right)$ in \eqref{cond-a1}, we obtain
\begin{equation}\label{cond-a2}
     -\rho_1 i  \xi_n \left(\Phi_n  \left(x- \frac{l}{2} \right) , \varphi_{n,x} + \psi_n\right)_{L^2} - \rho_1\left(\Phi_n  \left(x- \dfrac{l}{2} \right), \Phi_{n,x} + \Psi_n\right)_{L^2}\rightarrow 0.
\end{equation}
Computing \eqref{cond-a} + \eqref{cond-a2}, we get
  \begin{equation}\label{cond-a3}
    k\underbrace{\left((\varphi_{n,x} + \psi_n)_x  \left(x - \frac{l}{2}\right), \varphi_{n,x} + \psi_n \right)_{L^2}}_{:=R_n} +\,  \rho_1\left(\Phi_n  \left(x- \dfrac{l}{2} \right), \Phi_{n,x} + \Psi_n\right)_{L^2} \rightarrow 0.
\end{equation}
Integrating by parts, we obtain
\begin{eqnarray}\label{cond-a4}
\nonumber \re \; \left(\Phi_n  \left(x - \dfrac{l}{2}\right), \Phi_{n,x} \right)_{L^2} &=& \dfrac{1}{2}\int_{0}^{l}{\left(x -\frac{l}{2} \right)\frac{d}{dx} |\Phi_n|^2 dx}\\
%\nonumber &=& \dfrac{1}{2} \left(x - \frac{l}{2} \right)|\Phi_n|^2\bigg|_{0}^{l} - \dfrac{1}{2} \|\Phi_n\|_{L^2}^2\\
&=& \dfrac{l}{4} (|\Phi_n(0)|^2 + |\Phi_n(l)|^2) - \dfrac{1}{2} \|\Phi_n\|_{L^2}^2.
\end{eqnarray}
Similarly, we get
\begin{equation}\label{cond-a5}
    \re (R_n )= \dfrac{l}{4} (|\varphi_{n,x}(0) + \psi_n(0)|^2 + |\varphi_{n,x}(l) + \psi_n(l)|^2) - \dfrac{1}{2}\|\varphi_{n,x} + \psi_n\|_{L^2}^2.
\end{equation}
Taking the real part in \eqref{cond-a3} and using the identities given in \eqref{cond-a4}-\eqref{cond-a5}, we obtain
\begin{eqnarray}\label{cond-a6}
   \nonumber && \dfrac{\rho_1 l}{4} (|\Phi_n(0)|^2 + |\Phi_n(l)|^2) - \dfrac{\rho_1}{2} \|\Phi_n\|_{L^2}^2 + \dfrac{k l}{4} (|\varphi_{n,x}(0) + \psi_n(0)|^2 + |\varphi_{n,x}(l) + \psi_n(l)|^2)\\
    && \qquad \qquad \qquad \qquad \qquad \qquad \qquad \qquad \qquad \qquad \qquad - \dfrac{k}{2}\|\varphi_{n,x} + \psi_n\|_{L^2}^2 \rightarrow 0.
\end{eqnarray}

Combining \eqref{conv-norma} and Lemmas \ref{e-seq-tem-resultado1-compl}-\ref{e-seq-temp-resultado5}, we deduce
\begin{equation}\label{cond-a6-1}
    \dfrac{k}{2}\|\varphi_{n,x} + \psi_n\|_{L^2}^2 + \dfrac{\rho_1}{2}\|\Phi_n\|_{L^2}^2 \rightarrow \dfrac{1}{2}.
\end{equation}
    Computing \eqref{cond-a6}+\eqref{cond-a6-1}, we conclude
\begin{equation}\label{cond-a611}
\dfrac{l}{2} \left(\rho_1|\Phi_n(0)|^2 + \rho_1|\Phi_n(l)|^2  + k|\varphi_{n,x}(0) + \psi_n(0)|^2 + k|\varphi_{n,x}(l) + \psi_n(l)|^2\right)\rightarrow 1.
\end{equation}

We claim that $(\varphi_{n,x}(0))_{n\in \mathbb{N}}$ and $(\varphi_{n,x}(l))_{n\in \mathbb{N}}$ are bounded. Indeed, let us consider $p\in [0,l]$. Due to continuous embedding $H_0^1(0,l)\subset C([0,l])$, we infer $|\psi_n(p)|\leq C\|\psi_{n,x}\|_{L^2}$. Thus, from Lemma \ref{e-seq-temp-resultado3} we have $|\psi_n(0)|, \, |\psi_n(l)|\rightarrow 0$. Further, since $|\varphi_{n,x}(p)|^2\leq 2|\varphi_{n,x}(p)+ \psi_{n}(p)|^2 + 2|\psi_n(p)|^2$ in particular for $p=0$ or $p=l$, the desired follows from \eqref{cond-a611}.
Therefore, from convergences in \eqref{front-psi+int}, we obtain 
\begin{equation*}
    \left(\beta \psi_{n,x} + \int_{0}^{\infty}{g(s) \eta_{n, x}(s) ds}\right) \overline{\varphi_{n,x}}\bigg|_{0}^{l}\rightarrow 0.
\end{equation*}
Going back to \eqref{e-seq-temp-resultado-front} we obtain $\|\varphi_{n,x} + \psi_n\|_{L^2}^2 \rightarrow 0$  and, from Lemma \ref{e-seq-temp-resultado5}, we have   $\rho_1 \|\Phi_n\|_{L^2}^2 \rightarrow 0$. From this and again from  Lemmas \ref{e-seq-tem-resultado1-compl}-\ref{e-seq-temp-resultado5}, we finally concluded 
\begin{equation*}
    \|U_n\|^2_{\mathcal{H}_1} \rightarrow 0.
\end{equation*}
\qed

\subsubsection{Resolvent's upper bound}

In order to prove the second property in \eqref{properties}, it is enough to prove that 
\begin{equation}\label{bound-resol-eq}
\| U \|_{\mathcal{H}_j}=\| (i \lambda I_d - \mathcal{A}_j)^{-1} F \|_{\mathcal{H}_j} \leq C \|F\|_{\mathcal{H}_j}, \ \mbox{ as } \  |\lambda|\to\infty, \ j=1,2,
\end{equation}
for any $F\in \mathcal{H}_j$. Thus, given $F=(f_{1},f_{2},f_{3},f_{4},f_{5},f_{6}, f_{7})\in \mathcal{H}_j,$  let us start by considering  the resolvent equation 
\begin{equation} \label{resolvent-eq}
i\lambda U- {\cal A}_j U=F, \quad j=1,2,
\end{equation}
with $U=(\varphi, \Phi, \psi, \Psi,\theta, \Theta, \eta)\in D({\cal A}_j)$, where we recall that  ${\cal A}_j$ is defined in (\ref{semigroup 3.2}). In addition, equation \eqref{resolvent-eq} written  in terms of its components turns into
\begin{eqnarray}
i\lambda\varphi-\Phi & = & f_1 ,  \label{comp1} \\
i\lambda\rho_1\Phi-k(\varphi_x+\psi)_x +\sigma \,\Theta_x&=& \rho_1f_2 , \label{comp2} \\
i\lambda\psi-\Psi &=& f_3 , \label{comp3} \\
i\lambda\rho_2\Psi + k(\varphi_x+\psi) - \beta\psi_{xx} - \int_{0}^{\infty}{g(s)\eta_{xx}(s) ds} - \sigma \,\Theta & = & \rho_2f_4 ,\label{comp4} \\
i \lambda \theta - \Theta &=& f_5, \label{comp5}\\
i\lambda\rho_3\Theta-\delta\theta_{xx}- \gamma \Theta_{xx} +\sigma(\Phi_x+\Psi) & = & \rho_3f_6, \label{comp6}\\
 i\lambda \eta + \eta_s -\Psi & = & f_7. \label{comp7}
\end{eqnarray}

The proof of \eqref{bound-resol-eq} will be done as a consequence of some lemmas provided below in combination with the observability result, more precisely  Corollary \ref{cor-3}.
Hereafter, we  simplify the notation by using the same parameter $C>0$ to denote   different  constants. Besides, H\"{o}lder and Poincar\'e's inequalities, as well as $|\lambda|>1$ large enough, 
will be used   several times, possibly with no mention.

\begin{Lemma}\label{fisrt-estimate}
Under the assumptions of Theorem $\ref{theo-expo3}$ and the above notations,   there exists a constant $C>0$ such that
\begin{equation}\label{est-theta-x}
   \int_{0}^{\infty}{[- g'(s)] \|\eta_x(s)\|^2_{L^2}ds}, \, \|\Theta_{x}\|^2_{L^2} , \,   \|\eta\|_{\mathcal{M}}^2 \leq {C}\|U\|_{{\cal H}_j}\|F\|_{{\cal H}_j}, \quad j=1,2.
\end{equation}
\end{Lemma}
\begin{proof}
 Taking the inner product of $U$ with \eqref{resolvent-eq} in $\mathcal{H}_j$, and using \eqref{diss1}, we have
\begin{equation}\label{first-estimate-eq1}
\gamma \|\Theta_x\|_{L^2}^2 + \frac{1}{2} \int_{0}^{\infty}{[- g'(s)] \|\eta_x(s)\|^2_{L^2} ds} \, \leq \, \|U\|_{\mathcal{H}_j} \|F\|_{\mathcal{H}_j},
\end{equation}
from where it follows the first two estimates in \eqref{est-theta-x}. Also, from  \eqref{sobreg} and \eqref{first-estimate-eq1}, one easily concludes that
\begin{equation*}
\|\eta\|_{\mathcal{M}}^2\leq \dfrac{1}{k_1}\int_{0}^{\infty}{[-g'(s)]\|\eta_x(s)\|_{L^2}^2ds}\leq \dfrac{2}{k_1}\|U\|_{\mathcal{H}_j} \|F\|_{\mathcal{H}_j}.
\end{equation*}
\end{proof}

\begin{Lemma}\label{fisrt-estimate-2}
Under the assumptions of Theorem $\ref{theo-expo3}$ and the above notations,   there exists a constant $C>0$ such that 
\begin{equation}\label{est-theta-x-2}
    \|\theta_{x}\|^2_{L^2},\, \|\delta\theta_{x} + \gamma \Theta_x\|^2_{L^2}  \leq {C}\|U\|_{{\cal H}_j}\|F\|_{{\cal H}_j} , \quad j=1,2.
\end{equation}
\end{Lemma}
\begin{proof}
Deriving (\ref{comp5}) and taking the multiplier $\overline{\theta_x}$ in the resulting equation, we have
\begin{equation*}
i\lambda\int_0^l|\theta_x|^2\;dx=\int_0^l\Theta_x \overline{\theta_x}\;dx + \int_0^l f_{5,x} \overline{\theta_x}\;dx.
\end{equation*}
From Young's inequality and (\ref{est-theta-x}), we deduce the first estimate in (\ref{est-theta-x-2}), with $|\lambda|>1$ large enough. The estimate for $\delta \theta_x + \gamma \Theta_x$ follows from the first part and (\ref{est-theta-x}). 
\end{proof}

\begin{Lemma}\label{eLema3}
Under the assumptions of Theorem $\ref{theo-expo3}$ and the above notations,   there exists a constant $C>0$ such that
	\begin{equation}\label{eLema3-ii-complet}
	\rho_2 \|\Psi\|^2_{L^2}\leq C \|U\|_{\mathcal{H}_j} \|F\|_{\mathcal{H}_j} + C \|\eta\|_{\cal M}\left(\dfrac{\|\Phi\|_{L^2}}{|\lambda|}+ \|\psi_x\|_{L^2}\right), \quad j=1,2.
	\end{equation}
\end{Lemma}
\begin{proof}
	Taking the multiplier $\int_{0}^{\infty}{g(s) \overline{\eta(s)}ds}$ in \eqref{comp4} and performing integration by parts, we get
	\begin{eqnarray*}
		\rho_2 \int_{0}^{l}{\int_{0}^{\infty}{g(s) f_4 \overline{\eta(s)}ds}dx} &=& \underbrace{i\lambda \rho_2 \int_{0}^{l}{\int_{0}^{\infty}{g(s) \Psi\overline{\eta(s)}ds}dx}}_{:=R_{1}} - \sigma\int_{0}^{l}{\int_{0}^{\infty}{g(s) \Theta \overline{\eta(s)} ds}dx}  \\
		&& + \beta \int_{0}^{l}{\int_{0}^{\infty}{g(s) \overline{\eta_x(s)}\psi_xds}dx} +  \int_{0}^{l}{\left|\int_{0}^{\infty}{g(s)\eta_x(s) ds} \right|^2dx} \\
		&& + k \int_{0}^{l}{\int_{0}^{\infty}{g(s)(\varphi_x + \psi) \overline{\eta(s)}ds}dx}.
	\end{eqnarray*}
	Using \eqref{comp7} in $R_{1}$, we obtain
	\begin{eqnarray}\label{eLema3-1-compl}
	\nonumber &&\rho_2 \int_{0}^{l}{\int_{0}^{\infty}{g(s) |\Psi|^2ds}dx}=  \rho_2 \int_{0}^{l}{\int_{0}^{\infty}{g(s) \Psi \overline{f_5}ds}dx} - \underbrace{\rho_2 \int_{0}^{l}{\int_{0}^{\infty}{g(s) \Psi \overline{\eta_s(s)}ds}dx}}_{:=R_{2}}\\
	    && -\beta \int_{0}^{l}{\int_{0}^{\infty}{g(s) \overline{\eta_x(s)}\psi_xds}dx} - \underbrace{\int_{0}^{l}{\left|\int_{0}^{\infty}{g(s)\eta_x(s) ds} \right|^2dx}}_{:=R_{3}}\\
	\nonumber   &&  - \underbrace{k \int_{0}^{l}{\int_{0}^{\infty}{g(s)(\varphi_x + \psi) \overline{\eta(s)}ds}dx}}_{:=R_{4}} + \rho_2 \int_{0}^{l}{\int_{0}^{\infty}{g(s) f_4 \overline{\eta(s)}ds}dx}\\
	\nonumber &&+ \sigma\int_{0}^{l}{\int_{0}^{\infty}{g(s) \Theta \overline{\eta(s)} ds}dx}.
	\end{eqnarray}
	Performing integration by parts in $R_{2}$, applying Cauchy-Schwarz and Young's inequalities, we have
	\begin{eqnarray}\label{eLema3-R1}
	|R_{2}|\leq \rho_2 \sqrt{b_0} \|\Psi\|_{L^2} \left( -\int_{0}^{\infty}{g'(s) \|\eta(s)\|^2_{L^2}ds} \right)^{\frac{1}{2}},
	\end{eqnarray}
where we have used the notation   $b_0:= \int_{0}^{\infty}{g(s)ds}$ (and then $\beta= b-b_0$). 	Moreover, applying H\"{o}lder inequality, we get
	\begin{equation}\label{eLema3-3}
	|R_3|\leq \int_{0}^{l}{\left(\int_{0}^{\infty}{g(s)|\eta_x(s)| ds} \right)^2 dx}\leq  b_0 \|\eta\|^2_{\mathcal{M}}.
	\end{equation}
Also, equations \eqref{comp1} and \eqref{comp3} yield
	\begin{eqnarray*}
		\nonumber R_{4}&=& - \dfrac{i k}{\lambda} \int_{0}^{l}{\int_{0}^{\infty}{g(s) \overline{\eta_x(s)} ds} \Phi dx} - \dfrac{i k}{\lambda} \int_{0}^{l}{\int_{0}^{\infty}{g(s) \overline{\eta_x(s)} ds} f_1 dx} \\
		&& +\dfrac{i k}{\lambda} \int_{0}^{l}{\int_{0}^{\infty}{g(s) \overline{\eta(s)} ds} f_3 dx} + \dfrac{ik}{\lambda}\int_{0}^{l}{\int_{0}^{\infty}{g(s) \overline{\eta(s)} ds} \Psi dx}.
	\end{eqnarray*}
	Using Poincar\'e's inequality, we get
	\begin{equation}\label{eLema3-i1-compl}
  |R_{4}| \leq  C \|\eta\|_{\cal M} \|\Phi\|_{L^2} + C \|\eta\|_{\cal M} \|f_{1,x} + f_3\|_{L^2} + \dfrac{C}{|\lambda|} \|\eta\|_{\cal M} \|f_{3,x}\|_{L^2} + \dfrac{C}{|\lambda|} \|\eta\|_{\cal M} \|\Psi\|_{L^2}.
	\end{equation}
	Applying  H\"{o}lder's inequality in \eqref{eLema3-1-compl} and estimates \eqref{eLema3-R1}-\eqref{eLema3-i1-compl}, we can estimate
	\begin{eqnarray}\label{eLema3-i2-compl}
	\nonumber \rho_2 \|\Psi\|_{L^2}^2 &\leq& C\|\Psi\|_{L^2}\|f_{5}\|_{\mathcal{M}} + C \|\Psi\|^2_{L^2} + C \int_{0}^{\infty}{[ -g'(s)] \|\eta_x(s)\|^2_{L^2}ds} + C\|\eta\|^2_{\mathcal{M}} \\
	\nonumber && + C\|\psi_x\|_{L^2} \|\eta\|_{\mathcal{M}} + C \|\eta\|_{\mathcal{M}} \|f_4\|_{L^2} + \dfrac{C}{|\lambda|} \|\eta\|_{\cal M} \|\Phi\|_{L^2} + C\|\eta\|_{\mathcal{M}} \|\Theta\|_{L^2}\\
	\nonumber && + \dfrac{C}{|\lambda|} \|\eta\|_{\cal M} \|f_{1,x} + f_3\|_{L^2} + \dfrac{C}{|\lambda|} \|\eta\|_{\cal M}  \|f_{3,x}\|_{L^2} + \underbrace{\dfrac{C}{|\lambda|} \|\eta\|_{\cal M} \left(\|\Psi\|_{L^2}\right)}_{:=R_{5}}.
	\end{eqnarray}
	Applying Young's inequality in $R_{5}$ its follow that
	\begin{equation*}
	\nonumber \rho_2 \|\Psi\|_{L^2}^2 \leq C\|U\|_{\mathcal{H}_j} \|F\|_{\mathcal{H}_j} + C\|\eta\|_{\mathcal{M}}\|\psi_x\|_{L^2} + \dfrac{C}{|\lambda|} \|\eta\|_{\cal M} \|\Phi\|_{L^2} + C\|\eta\|_{\mathcal{M}} \|\Theta\|_{L^2}.
	\end{equation*}
From this and  Lemma \ref{fisrt-estimate}, we finally arrive at \eqref{eLema3-ii-complet} as desired.
\end{proof}
\begin{Lemma}\label{eLema4-complet}
Under the assumptions of Theorem $\ref{theo-expo3}$ and the above notations,   there exists a constant $C>0$ such that
	\begin{equation}\label{eLema-psix}
	\beta \|\psi_x\|_{L^2}^2\leq C \|U\|_{\mathcal{H}_j} \|F\|_{\mathcal{H}_j} + \dfrac{C}{|\lambda|}\|\varphi_x + \psi\|_{L^2}\|\Psi\|_{L^2} + C \|\Psi\|^2_{L^2}, \quad j=1,2.
	\end{equation}
\end{Lemma}
\begin{proof} 
	Taking the multiplier $\overline{\psi}$ in \eqref{comp4} we obtain
	\begin{eqnarray*}\label{eLema4-1-complet}
		\nonumber && \underbrace{i\lambda \rho_2 \int_{0}^{l}{\Psi \overline{\psi} dx}}_{:=S_{1}} +  \beta \int_{0}^{l}{ |\psi_x|^2 dx} + \int_{0}^{l}{\int_{0}^{\infty}{g(s) \eta_x(s) \overline{\psi_x}ds}dx} + k \int_{0}^{l}{\varphi_x \overline{\psi}dx}\\
		&& +  k\int_{0}^{l}{|\psi|^2dx} - \sigma \int_{0}^{l}{\Theta \overline{\psi} dx}= \rho_2 \int_{0}^{l}{f_4 \overline{\psi}dx}.
	\end{eqnarray*}
	Replacing  \eqref{comp3} in $S_{1}$, we get
	\begin{eqnarray*}\label{eLema4-2-complet}
		\nonumber  \beta \int_{0}^{l}{ |\psi_x|^2 dx} &=&  - \int_{0}^{l}{\int_{0}^{\infty}{g(s) \eta_x(s) \overline{\psi_x}ds}dx} - \underbrace{k \int_{0}^{l}{(\varphi_x + \psi) \overline{\psi}dx}}_{:=S_{2}}\\
		&& + \rho_2 \int_{0}^{l}{f_4 \overline{\psi}dx}  + \rho_2 \int_{0}^{l}{\Psi \overline{f_3}dx} + \rho_2\int_{0}^{l}{|\Psi|^2dx}\\
		\nonumber && +\sigma \int_{0}^{l}{\Theta \overline{\psi} dx}.
	\end{eqnarray*}
Inserting now \eqref{comp3} in $S_{2}$, we obtain
	\begin{eqnarray*}
		\nonumber  \beta \int_{0}^{l}{ |\psi_x|^2 dx} &=&  - \int_{0}^{l}{\int_{0}^{\infty}{g(s) \eta_x(s) \overline{\psi_x}ds}dx} + \dfrac{ik}{\lambda} \int_{0}^{l}{(\varphi_x + \psi) \overline{\Psi}dx} +\sigma \int_{0}^{l}{\Theta \overline{\psi} dx} \\
		&& + \dfrac{ik}{\lambda} \int_{0}^{l}{(\varphi_x + \psi) \overline{f_3}dx} + \rho_2 \int_{0}^{l}{f_4 \overline{\psi}dx}  + \rho_2 \int_{0}^{l}{\Psi \overline{f_3}dx} + \rho_2\int_{0}^{l}{|\Psi|^2dx}.
	\end{eqnarray*}
From the  above identity we can conclude \eqref{eLema-psix}  for some constant $C> 0$ and  $|\lambda|>1$ large enough.
\end{proof} 

\begin{Corollary}
	\label{cor-psi-concl}
Given $\varepsilon> 0$, there exists  $C_{\varepsilon}> 0$ such that
	\begin{equation}\label{eq-cor}
	\|\psi_x\|_{L^2}^2\leq \varepsilon \|U\|^2_{\mathcal{H}_j} + C_{\varepsilon}\|F\|^2_{\mathcal{H}_j}, \quad j=1,2.
	\end{equation}
\end{Corollary}
\begin{proof}
Combining Lemmas \ref{eLema3} and \ref{eLema4-complet}, we get
	\begin{eqnarray*}\label{psi-concl}
		\nonumber \|\psi_x\|_{L^2}^2 &\leq&  \dfrac{C}{|\lambda|} \| \varphi_x + \psi\|_{L^2} \left[\|F\|_{\mathcal{H}_j}^{\frac{1}{2}} \|U\|_{\mathcal{H}_j}^{\frac{1}{2}} + \|\eta\|_{\mathcal{M}}^{\frac{1}{2}} \left(\dfrac{\|\Phi\|_{L^2}^{\frac{1}{2}}}{|\lambda|^{\frac{1}{2}}} + \|\psi_x\|_{L^2}^{\frac{1}{2}} \right) \right]\\
		&& + C\|\eta\|_{\mathcal{M}} \left(\dfrac{\|\Phi\|_{L^2}}{|\lambda|} + \|\psi_x\|_{L^2}\right)+ C \|U\|_{\mathcal{H}_j} \|F\|_{\mathcal{H}_j}\\
	&\leq&  C \|U\|_{\mathcal{H}_j}^{\frac{3}{2}} \|F\|_{\mathcal{H}_j}^{\frac{1}{2}} + C\|U\|_{\mathcal{H}_j}^{\frac{3}{2}} \|\eta\|^{\frac{1}{2}}_{\mathcal{M}} + C \|U\|_{\mathcal{H}_j} \|\eta\|_{\mathcal{M}} + C \|U\|_{\mathcal{H}_j} \|F\|_{\mathcal{H}_j},
	\end{eqnarray*}
for some $C>0$. Therefore, from  Young's inequality and Lemma \ref{fisrt-estimate} along with Young's inequality once more, one can conclude \eqref{eq-cor}.
% 
%	\begin{equation}%\label{psi-concl2}
%	\|\psi_x\|_{L^2}^2\leq C \|U\|_{\mathcal{H}_j}^{\frac{3}{2}} \|F\|_{\mathcal{H}_j}^{\frac{1}{2}} + C\|U\|_{\mathcal{H}_j}^{\frac{3}{2}} \|\eta\|^{\frac{1}{2}}_{\mathcal{M}} + C \|U\|_{\mathcal{H}_j} \|\eta\|_{\mathcal{M}} + C \|U\|_{\mathcal{H}_j} \|F\|_{\mathcal{H}_j},
%	\end{equation} 
\end{proof}

In the next result, in order to keep the estimates for both cases $j=1$ and $j=2$, we are going to work with local estimates by means of cut-off functions to avoid point-wise terms. 
This procedure has been used e.g. in \cite{alves-marcio-matofu-rivera,michele-JEE} and here we adopt the same strategy.

\smallskip 
 
Let us consider  $l_0\in(0,l)$ and  $\delta_0>0$  arbitrary numbers     such that $(l_0-\delta_0,l_0+\delta_0)\subset(0,l),$ and a function $s\in C^2(0,l)$  satisfying 
\begin{equation}\label{hips1}
\mbox{supp }s\subset (l_0-\delta_0 , l_0+ \delta_0),  \quad 0\leq s(x) \leq 1, \ \ x\in(0,l),
\end{equation}
and
\begin{equation}\label{hips2}
s(x)=  1\quad\mbox{for}\quad x\in[l_0-\delta_0/2 , l_0+ \delta_0/2].
\end{equation}

\begin{Lemma}\label{Lphi}
Under the assumptions of Theorem $\ref{theo-expo3}$ and the above notations,   there exists a constant $C>0$ such that
\begin{eqnarray}\label{est-phix+Phi1}
\nonumber \int_{l_0-\frac{\delta_0}{2}}^{l_0 + \frac{\delta_0}{2}}{\left(|\varphi_x + \psi|^2 + |\Phi|^2\right)dx}&\leq& \dfrac{C}{|\lambda|^{\frac{3}{2}}}\left( \|\delta \theta_x + \gamma \Theta_x\|_{L^2}^{\frac{1}{2}}\|U\|_{\mathcal{H}_j}^{\frac{1}{2}} +  \|U\|_{\mathcal{H}_j}^{\frac{1}{2}}\|F\|_{\mathcal{H}_j}^{\frac{1}{2}}\right) \|U\|_{\mathcal{H}_j} \\
 && + \dfrac{C}{|\lambda|^{\frac{4}{3}}} \|\delta \theta_x + \gamma \Theta_x\|_{L^2}^{\frac{1}{2}}\|U\|_{\mathcal{H}_j}^{\frac{4}{3}} + C\|U\|_{\mathcal{H}_j} \|F\|_{\mathcal{H}_j} \\
\nonumber&& + \dfrac{C}{|\lambda|} \|F\|^2_{\mathcal{H}_j} + \dfrac{C}{|\lambda|} \|\delta \theta_x + \gamma \Theta_x\|_{L^2}\|U\|_{\mathcal{H}_j}, \quad j=1,2.
\end{eqnarray}
In addition,  given  $\epsilon>0$ there exists a constant $C_{\epsilon}>0$ such that
\begin{equation} \label{est-PHI-eps}
\int_{l_0-\frac{\delta_0}{2}}^{l_0+\frac{\delta_0}{2}}{(|\varphi_x + \psi|^2+|\Phi|^2)}\;dx \, \leq \,
\epsilon \|U\|_{{\cal H}_j}^2+C_\epsilon\|F\|_{{\cal H}_j}^2, \quad j=1,2.
\end{equation}
\end{Lemma}
\begin{proof}
Computing 
  $(\ref{comp1})_x + (\ref{comp3})$ and inserting the resulting expression in (\ref{comp6}), we get
\begin{equation}\label{comp5-modi2}
    i\lambda \rho_3 \Theta - (\delta \theta + \gamma \Theta)_{xx} + i\lambda \sigma (\varphi_x + \psi)= \rho_3 f_6 + \sigma (f_{1,x} + f_3).
\end{equation}
Taking the multiplier $sk\overline{[\varphi_x+\psi]}$ in 
(\ref{comp5-modi2}), integrating on $(0, l)$ and performing integration by parts, we have
\begin{align}\label{ob-lema31-conseq}
\nonumber & \underbrace{i\lambda \rho_3 \int_0^{l} \Theta ks\overline{(\varphi_x +\psi)} \;dx}_{:= I_1} - \underbrace{\int_0^{l} (\delta \theta_x + \gamma \Theta_x)_{x}ks\overline{(\varphi_x + \psi)} \;dx}_{:=I_2} \\
 &\quad =- i\lambda \sigma \int_0^{l} (\varphi_x + \psi) ks \overline{ (\varphi_x + \psi)} \;dx + \rho_3 \int_0^l f_6 ks\overline{(\varphi_x + \psi)}\; dx\\
\nonumber &\quad \quad + \sigma \int_0^l (f_{1,x} + f_3) ks \overline{(\varphi_x + \psi)}\; dx.
\end{align}
From equations (\ref{comp1}) and (\ref{comp3}), we achieve
\begin{equation*}
    I_1=\rho_3 \int_0^l (ks \Theta)_x \overline{\Phi}\; dx - \rho_3 \int_0^l \Theta ks\overline{\Psi}\; dx -\rho_3\int_0^l \Theta ks\overline{(f_{1,x} + f_3)}\; dx.
\end{equation*}
Performing again integration by parts, we have
\begin{eqnarray*}
    \nonumber I_2&=&  - \int_0^l (\delta \theta_x + \gamma \Theta_x) ks_x \overline{(\varphi_x + \psi)}\; dx \underbrace{- \int_0^l (\delta \theta_x + \gamma \Theta_x) ks \overline{(\varphi_x + \psi)_x}\; dx}_{:=I_3}.
\end{eqnarray*}
Using \eqref{comp2} one gets
\begin{equation*}
I_3=-\int_0^l (\delta \theta_x + \gamma \Theta_x) s\overline{(i\lambda \rho_1 \Phi)}\; dx - \int_0^l (\delta \theta_x + \gamma \Theta_x) \sigma s\overline{\Theta_x}\; dx + \int_0^l (\delta\theta_x + \gamma \Theta_x)\rho_1 s\overline{f_2}\; dx.
\end{equation*}
Replacing the identity provided by $I_3$ in $I_2$, and then replacing the identities provided in $I_2$ and $I_1$ in \eqref{ob-lema31-conseq}, we obtain
\begin{align*}
\nonumber i\lambda \int_0^l (\varphi_x + \psi) ks \overline{(\varphi_x + \psi)}\; dx=&\rho_3 \int_0^l f_6 ks\overline{(\varphi_x + \psi)}\; dx + \sigma \int_0^l (f_{1,x} + f_3) ks \overline{(\varphi_x + \psi)}\; dx\\
& - \int_0^l (\delta \theta_x + \gamma \Theta_x) ks_x\overline{(\varphi_x + \psi)}\; dx \\
& - \int_0^l (\delta \theta_x + \gamma \Theta_x) \sigma s \overline{\Theta_x}\; dx +\int_0^l (\delta \theta_x + \gamma \Theta_x) s\rho_1 \overline{f_2}\; dx\\
& - \rho_3 \int_0^l (ks\Theta)_x\overline{\Phi}\; dx + \rho_3\int_0^l \Theta ks\overline{\Psi}\; dx\\
& + \rho_3\int_0^l \Theta ks\overline{(f_{1,x} + f_3)}\; dx - \int_0^l (\delta \theta_x + \gamma \Theta_x) i\lambda s\rho_1 \overline{\Phi}\; dx.
\end{align*}
Going back to (\ref{ob-lema31-conseq}) and using  condition (\ref{hips1}) on $s$, we conclude
\begin{eqnarray}\label{ob-elema32}
|\lambda|\int_{l_0 - \delta_0}^{l_0 + \delta_0}{s |\varphi_x+\psi|^2 dx}&\leq& C|\lambda| \|\delta \theta_x + \gamma \Theta_x\|_{L^2} \left(\int_{l_0 - \delta_0}^{l_0 + \delta_0}{s |\Phi|^2 dx} \right)^{\frac{1}{2}} + C \|U\|_{\mathcal{H}_j} \|F\|_{\mathcal{H}_j} \\
\nonumber && + C\|\delta \theta_x + \gamma \Theta_x\|_{L^2} (\|F\|_{\mathcal{H}_j} +\|U\|_{\mathcal{H}_j}+\|\Theta_x\|_{L^2})+ C\|\Theta_x\|_{L^2} \|U\|_{\mathcal{H}_j}.\nonumber
\end{eqnarray}
Moreover, applying Young's  inequality and Lemmas \ref{fisrt-estimate}-\ref{fisrt-estimate-2} in \eqref{ob-elema32}, we obtain
\begin{align}
\int_{l_0-\delta_0}^{l_0+\delta_0}s|\varphi_x+\psi|^2\;dx \leq & \, C \, \|\delta \theta_x + \gamma \Theta_x\|_{L^2}\left(\int_{l_0-\delta_0}^{l_0+\delta_0}s|\Phi|^2 \, dx \right)^{1/2}  +  \frac{C}{|\lambda|}\|\delta \theta_x + \gamma \Theta_x\|_{L^2}\|U\|_{{\cal H}_j} \nonumber \\
&+ \frac{C}{|\lambda|}\|U\|_{{\cal H}_j}\|F\|_{{\cal H}_j} + \frac{C}{|\lambda|}\|F\|_{{\cal H}_j}^2 + \frac{C}{|\lambda|}\|\Theta_x\|_{L^2}\|U\|_{{\cal H}_j}. \label{est-phix+psi}
\end{align}

In what follows, we are going to estimate the term $\int_0^l{s}|\Phi|^2\;dx$. 
Indeed,  taking the multiplier $-{s}\overline{\varphi}$ in
(\ref{comp2}), performing integration by parts and applying (\ref{comp1}),  we get
\begin{align}\label{est-PHI1}
\rho_1\int_0^l{s}|\Phi|^2\;dx=k\int_0^l{s}|\varphi_x+\psi|^2\;dx- k \int_0^l {s}(\varphi_x+\psi) \overline{\psi}\;dx +I_4+I_5,
\end{align}
where
\begin{equation*}
I_4= 
\frac{i\sigma}{\lambda}\int_0^l s\, \Theta_x\overline{[\Phi+f_1]}\, dx -\rho_1\int_0^l {s}[\Phi  \overline{f_1}+f_2 \overline{\varphi}]\;dx \quad \mbox{and} \quad 
I_5=    k\int_0^l {s}'(\varphi_x+\psi) \overline{\varphi}\;dx .
\end{equation*}
It is easy to  see that
$$
|I_4| \leq \frac{C}{|\lambda|}\|\Theta_x\|_{L^2} \|U\|_{{\cal H}_j} + \frac{C}{|\lambda|}\|\Theta_x\|_{L^2}\|F\|_{{\cal H}_j}+ C \|U\|_{{\cal H}_j}\|F\|_{{\cal H}_j},
$$
for some constant $C>0$. In addition,   from equations (\ref{comp1}) and (\ref{comp3}), it follows that
$$
|\mbox{Re } I_5| \leq \frac{C}{|\lambda|}\|U\|_{{\cal H}_j}^2 + \frac{C}{|\lambda|}\|F\|_{{\cal H}_j}^2,
$$
for some constant $C>0$. Thus, taking the real part in (\ref{est-PHI1}) and  using \eqref{hips1},  we have
\begin{align*}\label{est-PHI2}
\int_{l_0-\delta_0}^{l_0+\delta_0}{s}|\Phi|^2\;dx \leq & \ C\int_{l_0-\delta_0}^{l_0+\delta_0} {s} |\varphi_x+\psi|^2\;dx +C\int_{l_0-\delta_0}^{l_0+\delta_0} s|\varphi_x+\psi| |\psi|\;dx  +\frac{C}{|\lambda|}\|\Theta_x\|_{L^2} \|U\|_{{\cal H}_j} \\
&   + \frac{C}{|\lambda|}\|\Theta_x\|_{L^2}\|F\|_{{\cal H}_j}+ C \|U\|_{{\cal H}_j}\|F\|_{{\cal H}_j} + \frac{C}{|\lambda|}\|U\|_{{\cal H}_j}^2 + \frac{C}{|\lambda|}\|F\|_{{\cal H}_j}^2 \\
\leq  & \ C\int_{l_0-\delta_0}^{l_0+\delta_0} {s} |\varphi_x+\psi|^2\;dx  +C  \left(\int_{l_0-\delta_0}^{l_0+\delta_0} {s} |\varphi_x+\psi|^2\;dx\right)^{1/2}\|\psi\|_{L^2}  
\\
& +  \frac{C}{|\lambda|}\|\Theta_x\|_{L^2}\|U\|_{{\cal H}_j} + \frac{C}{|\lambda|}\|\Theta_x\|_{L^2}\|F\|_{{\cal H}_j} + C\|U\|_{{\cal H}_j}\|F\|_{{\cal H}_j}  + \frac{C}{|\lambda|}\|U\|_{{\cal H}_j}^2\\
&  + \frac{C}{|\lambda|}\|F\|_{{\cal H}_j}^2.
\end{align*}
From estimate (\ref{est-phix+psi}), Young's  inequality, (\ref{est-theta-x}) and (\ref{est-theta-x-2}), it follows that
\begin{align*}
\int_{l_0-\delta_0}^{l_0+\delta_0}{s}|\Phi|^2\;dx \leq & \ C \, \|\delta \theta_x + \gamma \Theta_x\|_{L^2}\left(\int_{l_0-\delta_0}^{l_0+\delta_0}s|\Phi|^2 \, dx \right)^{1/2} + \frac{C}{|\lambda|}\|\delta \theta_x + \gamma \Theta_x\|_{L^2}\|U\|_{{\cal H}_j} \\
&  \,  +  \frac{C}{|\lambda|^{1/2}}\left(\|\delta \theta_x + \gamma \Theta_x\|_{L^2}^{1/2}\|U\|_{{\cal H}_j}^{1/2} +  \|U\|_{{\cal H}_j}^{1/2} \|F\|_{{\cal H}_j}^{1/2} +  \|F\|_{{\cal H}_j}\right)\|\psi\|_{L^2} \\
& +C\|F\|_{{\cal H}_j}^2  + \frac{C}{|\lambda|^2}\|U\|_{{\cal H}_j}^2  +  C \, \|\delta \theta_x + \gamma \Theta_x\|_{L^2}^{1/2}\left(\int_{l_0-\delta_0}^{l_0+\delta_0}s|\Phi|^2 \, dx \right)^{1/4}\|\psi\|_{L^2}\\
& + \frac{C}{|\lambda|^{1/2}} \|\Theta_x\|_{L^2}^{1/2} \|U\|_{{\cal H}_j}^{1/2}\|\psi\|_{L^2} + \frac{C}{|\lambda|}\|\Theta_x\|_{L^2} \|U\|_{{\cal H}_j} + \frac{C}{|\lambda|}\|\Theta_x\|_{L^2} \|F\|_{{\cal H}_j}\\
&  + C\|U\|_{{\cal H}_j}\|F\|_{{\cal H}_j}.
\end{align*}
Using again Young  inequality, (\ref{est-theta-x}), (\ref{est-theta-x-2}), and then equation  (\ref{comp3}), we get
\begin{align*}
\int_{l_0-\delta_0}^{l_0+\delta_0}{s}|\Phi|^2\;dx \leq & \,  \frac{C}{|\lambda|^{3/2}}\left(\|\delta \theta_x + \gamma \Theta_x\|_{L^2}^{1/2}\|U\|_{{\cal H}_j}^{1/2} +  \|U\|_{{\cal H}_j}^{1/2} \|F\|_{{\cal H}_j}^{1/2} +  \|F\|_{{\cal H}_j}\right)\big(\|U\|_{{\cal H}_j}+\|F\|_{{\cal H}_j}\big)\\ & \, + \frac{C}{|\lambda|^{4/3}} \|\delta \theta_x + \gamma \Theta_x\|_{L^2}^{2/3}\big(\|U\|_{{\cal H}_j}^{4/3}+\|F\|_{{\cal H}_j}^{4/3}\big)  + C\|U\|_{{\cal H}_j}\|F\|_{{\cal H}_j} \\
&+C\|F\|_{{\cal H}_j}^2  + \frac{C}{|\lambda|}\|U\|_{{\cal H}_j}^2.
\end{align*}
Applying once more Young's  inequality and estimate (\ref{est-theta-x-2}), we obtain
\begin{eqnarray}
\int_{l_0-\delta_0}^{l_0+\delta_0}s|\Phi|^2\;dx \, &\leq  & \ \frac{C}{|\lambda|^{3/2}} \left( \|\delta \theta_x + \gamma \Theta_x\|_{L^2}^{1/2}\|U\|_{{\cal H}_j}^{1/2}+ \|U\|_{{\cal H}_j}^{1/2}\|F\|_{{\cal H}_j}^{1/2} \right)\|U\|_{{\cal H}_j} \nonumber \\
&& + \frac{C}{|\lambda|^{4/3}} \|\delta \theta_x + \gamma \Theta_x\|_{L^2}^{2/3} \|U\|_{{\cal H}_j}^{4/3} + \frac{C}{|\lambda|}\|\delta \theta_x + \gamma \Theta_x\|_{L^2}\|U\|_{{\cal H}_j} \label{est-PHI} \\
&&      + C\|U\|_{{\cal H}_j} \|F\|_{{\cal H}_j} + C\|F\|_{{\cal H}_j}^2 + \frac{C}{|\lambda|}\|U\|_{{\cal H}_j}^2. \nonumber
\end{eqnarray} 
Therefore, combining  (\ref{est-phix+psi}) with (\ref{est-PHI}), using again Young's  inequality and estimate (\ref{est-theta-x-2}), we conclude
\begin{align*}
\int_{l_0-\delta_0}^{l_0+\delta_0}s\left(|\varphi_x+\psi|^2 + |\Phi|^2\right)dx \leq & \, \frac{C}{|\lambda|^{3/2}} \left( \|\delta \theta_x + \gamma \Theta_x\|_{L^2}^{1/2}\|U\|_{{\cal H}_j}^{1/2}+ \|U\|_{{\cal H}_j}^{1/2}\|F\|_{{\cal H}_j}^{1/2} \right)\|U\|_{{\cal H}_j}  \\
&   + \frac{C}{|\lambda|^{4/3}} \|\delta \theta_x + \gamma \Theta_x\|_{L^2}^{2/3}\|U\|_{{\cal H}_j}^{4/3}  + C\|U\|_{{\cal H}_j}\|F\|_{{\cal H}_j} + C\|F\|_{{\cal H}_j}^2,
\end{align*}
from where (\ref{est-phix+Phi1})-(\ref{est-PHI-eps}) follow by noting the 
properties of $s$ in (\ref{hips1})-(\ref{hips2}), Young's  inequality, and (\ref{est-theta-x-2}).
\end{proof}

We are now ready to prove \eqref{bound-resol-eq} as follows.

\smallskip 

\noindent \underline{\it Completion of the proof of \eqref{bound-resol-eq} for  $j=1,2.$}
Let  $\varepsilon> 0$ be given. 
Combining the estimates \eqref{eq-cor} and \eqref{est-PHI-eps} there exists a constant $C_{\varepsilon}$ such that
\begin{equation}\label{comp-proofs-10}
\int_{l_0-\frac{\delta_0}{2}}^{l_0+ \frac{\delta_0}{2}}  \big(|\varphi_x |^2+ |\Phi|^2\big)\;dx \, \leq \, \varepsilon \|U\|_{\mathcal{H}_j}^2+C_\varepsilon\|F\|_{\mathcal{H}_j}^2.
\end{equation}
This is the precise moment where we can take advantage of  Corollary \ref{cor-3} to extend the local estimate \eqref{comp-proofs-10} to the whole interval $(0,l)$. Indeed, through the first two components \eqref{comp1}-\eqref{comp2} of the resolvent equation, we observe that 
$V:= (\varphi, \Phi)$ is a solution of \eqref{eqondas1}-\eqref{eqondas2} with
$$g_1:=f_1 \quad \quad \mbox{and}\quad \quad g_2 :=\rho_1f_2-(\sigma\, \Theta)_x + k\psi_x,$$
and \eqref{comp-proofs-10} implies that (\ref{ass-corol}) is verified with $$b_1:=l_0- \frac{\delta_0}{2},  \ \    b_2:=l_0 + \frac{\delta_0}{2}, \ \  \mbox{and} \ \ \Lambda:=\varepsilon \|U\|_{\mathcal{H}_j}^2+C_\varepsilon\|F\|_{\mathcal{H}_j}^2 .$$ 
Thus, Corollary \ref{cor-3}, Lemma \ref{fisrt-estimate}, Corollary \ref{cor-psi-concl}, and Young's inequality imply
	\begin{equation}\label{last-exp}
	\int_{0}^{l}  \big(|\varphi_x|^2+ |\Phi|^2 \big)\;dx \, \leq \varepsilon C\|U\|_{\mathcal{H}_j}^2+C_\varepsilon \|F\|_{\mathcal{H}_j}^2 + C \|\psi_x\|_{L^2}^2 + C \|F\|^2_{\mathcal{H}_j},
	\end{equation}
	for some constants $C, C_\varepsilon>0$. 
	 Using Young and Poincar\'e's inequalities again, it is easy to see
\begin{equation}\label{comp-proofs-11}
\|\varphi_x + \psi\|^2_{L^2} \leq 2 \|\varphi_x\|_{L^2}^2 + 2 L^2 \|\psi_x\|^2_{L^2},
\end{equation}
and combining \eqref{last-exp}-\eqref{comp-proofs-11}, we have
\begin{equation}\label{last-exp2}
	\int_{0}^{l}  \big(|\varphi_x + \psi|^2+ |\Phi|^2 \big)\;dx \, \leq \varepsilon \|U\|_{\mathcal{H}_j}^2+C_\varepsilon \|F\|_{\mathcal{H}_j}^2 + C \|\psi_x\|_{L^2}^2,
	\end{equation}
	for some constants $C, C_\varepsilon>0$. Adding the estimates from \eqref{last-exp2}, and Lemmas \ref{eLema3} and \ref{eLema4-complet}, and using once again Young's inequality, we get	
\begin{eqnarray}\label{last-exp3}
	\int_{0}^{l}  \big(|\varphi_x + \psi|^2+ \beta |\psi_x|^2 + |\Phi|^2 + \rho_2 |\Psi|^2\big)\;dx &\leq& \varepsilon \|U\|_{\mathcal{H}_j}^2+C_{\varepsilon} \|F\|_{\mathcal{H}_j}^2 + C \|\psi_x\|_{L^2}^2 \\
\nonumber &&  + C \|F\|^2_{\mathcal{H}_j} + \dfrac{\varepsilon}{|\lambda|^2}\|\Phi\|_{L^2}^2 + C_{\varepsilon} \|\eta\|^2_{\mathcal{M}}.
	\end{eqnarray}
Finally, from  \eqref{last-exp3},   Lemmas \ref{fisrt-estimate} and \ref{fisrt-estimate-2}, and Corollary \ref{cor-psi-concl},  we achieve
	\begin{equation}\label{last-exp4}
\nonumber \|U\|_{\mathcal{H}_j}^2 \leq \varepsilon \|U\|_{\mathcal{H}_j}^2+C_{\varepsilon} \|F\|^2_{\mathcal{H}_j} + \dfrac{\varepsilon}{|\lambda|^2}\|\Phi\|_{L^2}^2.
	\end{equation}
Therefore, taking $\varepsilon>0$ small enough and $|\lambda|>1$ sufficiently large, we conclude that \eqref{bound-resol-eq} holds true.
This finishes the proof of Theorem \ref{theo-expo3}.
\qed

\begin{rem}\label{rem-santos}\rm
{\bf (a)} In the authors' opinion, 	Theorem \ref{theo-expo3} seems to correct the insight  claimed in   \cite[Sect. 7]{dilberto-mauro} with respect to 
the stability of a problem
 related to  system \eqref{intro 1.5}-\eqref{bc}, where it is stated
 (with no computations therein)
  that it is  not exponentially stable in general. Indeed, on p. 670 the authors claim {\it ``the model has the optimal polynomial decay rate when $\frac{\rho_{1}}{k} \neq \frac{\rho_{2}}{b},$ and that it is of the form $t^{-\frac{1}{2}}.$ Again, computations are required''.} 
	In the latter statement, 
	the authors refer to problem $ (7.1) $-$ (7.3) $ therein, which in turn corresponds to \eqref{intro 1.5}-\eqref{bc} unless mixed boundary conditions. Nevertheless, through the technical results employed in the proofs of the present paper,  Theorem \ref{theo-expo3} shows a different (and new) perspective to this assertion, namely, the  thermo-viscoelastic system \eqref{intro 1.5}-\eqref{bc} is exponentially stable,  and such uniform stability is independent of any relationship among the coefficients and boundary conditions addressed.
	
{\bf (b)} The physical reason for such uniform stability is already  highlighted in Remarks \ref{rem-strength}-\ref{rem-strength-a}.
Below we provide a  diagram that  clarifies the state of the art in the propagation of dissipativity along the solution according to the proofs of Lemmas \ref{fisrt-estimate} - \ref{Lphi}. It illustrates the strength of the thermo-(visco-)elastic damping feedback in problem \eqref{intro 1.5}-\eqref{bc}, by stressing the mathematical viewpoint of the damping propagation.

\begin{table}[h]
	\centering
	\begin{tabular}{c}
		{\bf  State of the art of damping propagation}  \medskip  \\
		\hline    \\                  
		The damping propagation over the bending moment is provided by the history \medskip \\
		$$
		\boxed{ \  \ \  {\color{blue}\eta} \ \ \ \stackrel{1^{st}}{\longmapsto} \ \ \ {\color{blue}  \Psi:=\psi_t} \ \ \ \stackrel{2^{nd}}{\longmapsto} \ \ \ {\color{blue} M=b\psi_{x}} \ \ \ }
	$$ \medskip \\
According to Lemmas \ref{fisrt-estimate},  \ref{eLema3}, and \ref{eLema4-complet}. 	\medskip 
	  \\ \hline \\
		The dissipation propagates to the shear force  through the thermal component
		  \medskip \\
	$$
	\boxed{
	\  \ \ {\color{red} \Theta} \ \ \ \stackrel{1^{st} }{\longmapsto} \ \ \ {\color{red} \theta}  \ \ \ \stackrel{2^{nd} }{\longmapsto} \ \ \ {\color{red} \Phi:=\varphi_t} \ \ \
	\stackrel{3^{rd}}{\longmapsto}	
	\ \ \ {\color{red} S=k(\varphi_{x}+\psi)}  \ \ \ 
	}
	$$	 \medskip \\
	According to Lemmas \ref{fisrt-estimate}, \ref{fisrt-estimate-2}, and \ref{Lphi}.
	\medskip 
	\\ \hline	
	\end{tabular}
\end{table}

\end{rem}

\appendix

\section{Appendix: auxiliary results}\label{sec-append-aux-results}

\subsection{Spectral results: a short review for linear operators}

In order to make this work more self-contained as possible  as well as for the reader's convenience, we report some  well-known technical   results on linear operators within the   functional analysis.

 We start with the following result that characterizes the spectrum $\sigma({A})=\mathbb{C} \backslash \rho({A})$ 
 of  linear operators $A:D(A)\subset X\to X$ whose domain $D(A)$ is compactly embedded in a Banach space $X.$  It can be found in the book by   Engel--Nagel \cite{Engel}.

 \begin{prop}[{\cite[Proposition 5.8 and Corollary 1.15]{Engel}}]\label{teo-inc-compacta} Let  $(X,\|\cdot\|_X)$ be a 
 	Banach space and  consider $A:D(A)\subset X\rightarrow X$ a linear operator with nonempty resolvent set   $\rho(A)$. 
 	\begin{enumerate}
 		\item[{\bf (i)}]   $A$ has compact resolvent (that is, there exists $\lambda\in\rho(A)$ such that $(\lambda I_d - A)^{-1}$ is compact) if  and only if  
 		the 	canonical injection
 		$i:(D(A),\|\cdot\|_{D(A)})\rightarrow (X,\|\cdot\|_X)$ is compact.
 		
 		\item[{\bf (ii)}] If operator  $A$ has compact resolvent, then the spectrum
 		$\sigma(A)$ consists only
 		of  eigenvalues of $A$.
 	\end{enumerate}	
 \end{prop}
 
As stated above, Proposition \ref{teo-inc-compacta} constitutes a very efficient tool when dealing with linear operators whose domain is compactly embedded in the space. The next result  illustrates this fact, which  has been hugely used in the literature concerning the stabilization of linear $C_0$-semigroups related to evolution models,    as one can see e.g.  in \cite{jaime-dilberto-santos-2014-1,alves-marcio-matofu-rivera,hugo-racke,Zhuangyi.Liu.book,santos-jde2012}.

\begin{Corollary}\label{theo-eixo-imag}
	Let $H$ be a Hilbert space and  $A: D(A)\subset H \rightarrow H$ be a linear operator with $\rho(A)\neq\emptyset$ and such that $i \lambda I_d-A$ is injective for every $\lambda\in\mathbb{R}$. If the embedding $D(A)\hookrightarrow H$ is compact, then   $i \mathbb{R}\subset \rho(A)$.
\end{Corollary}

\begin{proof}
Immediately from  Proposition \ref{teo-inc-compacta}.
\end{proof}  
 
  On the other hand, when 
  $D(A)$ is not (necessarily) compactly embedded in $H$, the picture changes considerable since this property is fundamental in the proof of   Proposition \ref{theo-eixo-imag}. Indeed,  the lack of compactness leads us to a more delicate way in the proof of the property $	i \mathbb{R} \, \subset \, \rho(A)$, as one can check in the proofs presented by \cite{conti,hugo-racke,Zhuangyi.Liu.book,jaime-hug0-jmaa08}. This second scenario is very common when one deals with viscoelastic problems driven by memory terms in the history context, where the domain $D(A)$ involves weighted spaces in terms of the memory component. See, for instance, \cite{giorgi-grasseli-pata99,grasseli-pata02,pata-zucchi} where is remarked that the embedding $D(A)\hookrightarrow H$ is not compact in general.
  
 To address this case on the non-compactness assumption, we rely on  Liu and Zheng's book \cite{Zhuangyi.Liu.book}. The next result can be stated as a consequence of the statements in  \cite[Chapt. 2]{Zhuangyi.Liu.book}.
 
\begin{prop}[{\cite[Sect. 2.2,  p. 25]{Zhuangyi.Liu.book}}]
\label{result-principal-apA}
	Let $H$ be a Hilbert space and  $A: D(A)\subset H \longrightarrow H$ be a linear closed operator
	with $\overline{D(A)}=H$ and such that $0\in \rho(A)$. Let us also set $\ell= \frac{1}{\|(-A)^{-1}\|_{\mathcal{L}(H)}}$. If $i\mathbb{R}\not \subset \rho(A)$, then there exist a real number $\omega \in (0, \ell]$, a sequence $\lambda_n\in \mathbb{R}$, with $|\lambda_n|<
	\omega$ and $|\lambda_n|\rightarrow \omega$, and a sequence $U_n\in D(A)$, with $\|U_n\|_{H}=1$, such that
	\begin{equation*}\label{Ap-resolve-conv-zero}
		(i\lambda_n I_d - A) U_n  \  \to \ 0 \ \ \mbox{ in } \ \ H.
	\end{equation*}
\end{prop} 

As an immediate  consequence (converse) of Proposition \ref{result-principal-apA}, we can write down. 
  
\begin{Corollary} \label{result-princ-apA-contrap}
	Let $H$	be a Hilbert space and  $A: D(A)\subset H \longrightarrow H$ be a linear closed operator
	with $\overline{D(A)}=H$ and such that $0\in \rho(A)$. Let us also set $\ell= \frac{1}{\|(-A)^{-1}\|_{\mathcal{L}(H)}}$. If, for every real number $\omega \in (0, \ell]$, every sequence $\lambda_n\in\mathbb{R}$, with $|\lambda_n|<
	\omega$ and $|\lambda_n| \omega$, and every sequence $U_n\in D(A)$, with $\|U_n\|_{H}=1$, we assume that
	\begin{equation*}
		%\label{Ap-resolve-conv-zero-contrap}
		(i\lambda_n I_d - A) U_n  \  \not\to \ 0 \ \ \mbox{ in } \ \ H,
	\end{equation*}
	then $i\mathbb{R} \subset \rho(A)$.
\end{Corollary}

\subsection{Observability inequality: a reading from the resolvent viewpoint}\label{sec-obser-ineq}

As  far as we know, the internal and boundary observability results for wave models are known   since 
the classical  book by Lions  \cite{lion-88}. Several other related results can be also found in Komornik \cite{komornik-94} employing the multiplier technique through the corresponding energy. Nowadays, there are several  well-established observability results related to wave systems
and their interplay between control and stability theory, and here we intend to reread such observability inequality suited to our case. In what concerns linear models involving wave-like systems, say Timoshenko models, we refer to \cite{alves-marcio-matofu-rivera,cacalvetal-zamp} where observability inequalities are provided. Thus, the following statements can be seen as a particular case of the results on observability presented by \cite{alves-marcio-matofu-rivera,cacalvetal-zamp,komornik-94,lion-88} from the resolvent equation viewpoint, and we review them here for didactic purposes.

Let us consider the following initial-boundary value problem related to the well-known one-dimensional  linear wave model
\begin{equation}\label{wave}
	\left\{\begin{array}{l}
		\varrho  u_{tt} - k u_{xx} = 0 \quad \mbox{in}\quad (0,L) \times (0,\infty), \medskip \\
		u(0,t)=u(L,t)=0,  \quad t\geq0, \medskip \\
		u(x,0)=u_0(x),  \ \  u_{t}(x,0)=u_1(x),
	\end{array}
	\right.
\end{equation}
where $u=u(x,t)$ represents the displacement of a vibrating string with length $L>0,$ $\varrho>0$ is the mass density, and $k>0$ stands for Young's modulus of the material. Thus,  $\sqrt{\frac {k}{\varrho }}$  stands for the 
speed of wave propagation  in the string. 

Model \eqref{wave} can be rewritten as in the first order abstract problem
\begin{equation} \label{abstract-wave}
	\left \{
	\begin{array}{l}
		V_{t}=\mathcal{B} V, \quad  t>0, \medskip \\
		V(0)=V_{0},
	\end{array}
	\right.
\end{equation}
where
$$ v:=u_t, \quad 
V:=\left[\begin{array}{c}
	u\\
	v 
\end{array}\right], \quad  
\mathcal{B} V=\left[\begin{array}{c}
	v\\
	\frac {k}{\varrho }u_{xx}
\end{array}\right], \quad V_0 = \left[\begin{array}{ccc}
	u_0\\
	u_1
\end{array}\right].
$$
We also set the Hilbert space 
$$
H=H^1_0(0,L)\times L^2(0,L)
$$
with standard norm 
$$
\|(u,v)\|_H^2 = \|u_x\|^2_{L^2} +\|v\|^2_{L^2}, \quad \forall \, (u,v)\in H,
$$
where $\|\cdot\|_{L^2}$ stands for the usual $L^2$-norm in $(0,L)$. In this case, the domain of operator $\mathcal{B}$ is given by
$$
D(\mathcal{B})=\left(H^2(0,L)\times H^1_0(0,L)\right)\times H^1_0(0,L).
$$
The above construction as well as 
the well-posedness of problem \eqref{abstract-wave} are very well-known in the literature, see for instance the classical books \cite{Zhuangyi.Liu.book,pazy}.

\smallskip

In what follows, we are going to proceed with the  observability and extension results by means of the resolvent equation corresponding to \eqref{abstract-wave}. To do so, 
let us  consider the resolvent equation 
\begin{equation}\label{resolvent-eq-wm}
	i\lambda V - \mathcal{B}V =G,
\end{equation}
for $\lambda\in \mathbb{R}$ and $G=(g_1,g_2)\in H$, which in terms of 
its components can be written as
\begin{eqnarray}
	i\lambda u-v= g_1\in H_0^1(0, L),&\label{eqondas1}  \\
	i\lambda \varrho v- k u_{xx}=\varrho g_2\in L^2(0, L).& \label{eqondas2}
\end{eqnarray}

 It is easy to show that  \eqref{resolvent-eq-wm} has a unique solution $V\in D(\mathcal{B})$ (cf. \cite{Zhuangyi.Liu.book}). 

 In what follows, given any numbers $0 \le a_1 < a_2 \le L$, the notation $\|V\|_{a_1,\,a_2}^2$ stands for
\begin{equation}\label{notation}
	\|V\|_{a_1,\,a_2}^2 =  \int_{a_1}^{a_2}\left(|u_x(x)|^2+|v(x)|^2 \right) dx \quad \mbox{ with } \quad \|V\|_{0,L}^2=\|V\|_{H}^2.
\end{equation}

\begin{prop}[Observability Inequality] \label{lema-observ}
	Under the above notations and taking $G=(g_1,g_2)\in H$, let $V=(u,v)$ be the (regular) solution of \eqref{resolvent-eq-wm}. If we consider any numbers  $0 \le a_1 < a_2 \le L$, then 
	there exist constants    $C_0, C_1>0$ (depending only on $\varrho$ and $k$) such that   
	\begin{equation} \label{wav1}
		|u_x(a_j)|^2 + |v(a_j)|^2  \leq C_0 \|V\|_{a_1, a_2}^2 + C_0 \|G\|_H^2 , \quad j=1,2,
	\end{equation}
	and
	\begin{align}\label{wav2}
		\|V\|_{a_1, a_2}^2 \leq C_1 \left[|u_x(a_j)|^2+|v(a_j)|^2\right] + C_1 \|G\|_H^2, \quad j=1,2.
	\end{align}
\end{prop} 
\begin{proof}	
It can be done as a particular case of 	\cite{alves-marcio-matofu-rivera,cacalvetal-zamp,komornik-94,lion-88}.	
\end{proof}

Two helpful extension results are given below  as a consequence of Proposition \ref{lema-observ}, which in turn will be the results used later in order to recover global estimates in the applications. Given  $b_1,b_2 \in [0,L]$, with $b_1 < b_2$, let us keep the notation \eqref{notation} in mind.

\begin{Corollary}\label{cor-3}
	Under the conditions of Proposition $ \ref{lema-observ} $, 
	let $V=(u,v)$ be the regular solution of \eqref{resolvent-eq-wm}.
	If for some sub-interval $(b_1,b_2) \subset (0,L)$ we have 
	\begin{equation}\label{ass-corol}
		\|V\|_{b_1,b_2}^2  \leq \Lambda, \quad \mbox{for some parameter } \Lambda=\Lambda(V,G,\lambda),
	\end{equation}
	then there exist a (universal) constant $C>0$ such that
	\begin{equation*} \label{concl-cor}
		\|V\|_{H}^2  \leq C \Lambda + C \|G\|_{H}^2.
	\end{equation*}
\end{Corollary}
\begin{proof}
	The proof is just a simple combination of (\ref{wav1})-(\ref{wav2}) along with hypothesis (\ref{ass-corol}).
	Indeed,  from (\ref{wav1}) and (\ref{ass-corol})  we have
	\begin{equation}  \label{prov-1}
		|u_x(b_j)|^2 + |v(b_j)|^2 \leq C_0 \Lambda + C_0 \|G\|_{0,\,L}^2, \quad j=1,2.
	\end{equation}
	On the one hand, using  (\ref{wav2}) with $a_1=0$, $a_2=b_2$ and (\ref{prov-1}) with $j=2$, we obtain
	$$
	\int_0^{b_2}\Big(|u_x|^2+|v|^2\Big)\,dx \leq C_2 \,\Lambda + C_2\|G\|_{H}^2,
	$$
	where $C_2 = C_1C_0+C_1>0$.  	
	On the other hand, using (\ref{wav2}) with $a_1 = b_2$, $a_2 = L$ and (\ref{prov-1}) with $j=1$,
	we also get
	$$
	\int_{b_2}^L\Big(|u_x|^2+|v|^2\Big)\,dx \leq C_2 \, \Lambda + C_2\|G\|_{H}^2.
	$$
	Therefore, adding these last two inequalities,  we arrive at the estimate (\ref{concl-cor}).
\end{proof}

\begin{Corollary}\label{cor-4}
	Under the conditions of Proposition $ \ref{lema-observ} $, 
	let $V=(u,v)$ be the regular solution of \eqref{resolvent-eq-wm}.
	If for some subinterval $(b_1,b_2) \subset (0,L)$ we have that
	\begin{equation}\label{ass-corol4}
		|u_x(b_j)|^2 + |v(b_j)|^2  \leq \Lambda, \quad j=1,2,
	\end{equation}
	then there exist a (universal) constant $C>0$ such that
	\begin{equation*} \label{concl-cor4}
		\|V\|_{H}^2  \leq C \Lambda + C \|G\|_{H}^2.
	\end{equation*}
\end{Corollary}
\begin{proof}
	Analogous to the proof of Corollary \ref{cor-3}, by using only (\ref{wav2}) and  (\ref{ass-corol4}).
\end{proof}

\paragraph{Acknowledgments.} The authors would like to express their gratitude to the anonymous referee for all remarks on a previous version as well as to the handling editor that allowed reconstructing this paper under a new perspective and improve itself with respect to the applied part.

\end{document}